\documentclass[a4paper]{article}

\usepackage[utf8x]{inputenc}
\usepackage{amsmath, amsfonts}

\usepackage[a4paper,top=3cm,bottom=2cm,left=3cm,right=3cm,marginparwidth=1.75cm]{geometry}

\usepackage{amsmath} 
\usepackage{graphicx} 
\usepackage{pdfpages}
\usepackage[colorlinks=true, allcolors=blue]{hyperref} 
\usepackage{xcolor}
\usepackage{authblk}
\usepackage{enumitem}
\usepackage{ulem}

\graphicspath{ {./figures/} }

\newtheorem{thm}{Theorem}[section]
\newtheorem{lemma}[thm]{Lemma}
\newtheorem{prop}[thm]{Proposition}

\newtheorem{rem}[thm]{Remark}

\newcommand{\R}{\ensuremath{\mathbb{R}}}
\def\C{{\mathbb{C}}}

\def\C{{\mathbb{C}}}

\def\sqr{{\rm sqrt}}
\def\Re{{\textrm{Re}}}
\def\Im{{\textrm{Im}}}
\def\cD{{\cal{D}}}

\numberwithin{equation}{section}
\newenvironment{pf}{\vspace{0.5em}\noindent\textbf{Proof.} }{\quad \hfill $\Box$ \\ \vspace{0.2em}\\}

\title{Long-time behaviour of the correlated random walk system}
\author[1]{Joaqu\'{\i}n Menacho}
\author[2]{Marta Pellicer}
\author[3]{J. Sol\`a-Morales}
\affil[1]{IQS School of Engineering, Universitat Ramon Llull, Via Augusta, 390, 08017 Barcelona, Catalunya, Spain, joaquin.menacho@iqs.edu}
\affil[2]{Departament d'Inform\`atica, Matem\`atica Aplicada i Estad\'istica, Universitat de Girona, EPS-4, c. Universitat de Girona, 6, 17003, Girona, Catalunya, Spain, marta.pellicer@udg.edu}
\affil[3]{Departament de Matem\`atiques, Universitat Polit\`ecnica de Catalunya, Av. Diagonal
647, 08028 Barcelona, Catalunya, Spain, jc.sola-morales@upc.edu}

\date{}
\begin{document}

\maketitle
\begin{abstract}
In this work, we consider the so-called correlated random walk system (also known as correlated motion or persistent motion system), used in biological modelling, among other fields, such as chromatography. This is a linear system which can also be seen as a weakly damped wave equation with certain boundary conditions. We are interested in the long-time behaviour of its solutions. To be precise, we will prove that the decay of the solutions to this problem is of exponential form, where the optimal decay rate exponent is given by the dominant eigenvalue of the corresponding operator. This eigenvalue can be obtained as a particular solution of a system of transcendental equations. A complete description of the spectrum of the operator is provided, together with a comprehensive analysis of the corresponding eigenfunctions and their geometry.
\end{abstract}

\section{Introduction and main results} \label{sec:intro}

The so-called correlated random walk system, sometimes known as the Goldstein-Kac system (see \cite{Goldstein,Kac}), has been widely used to model several biological processes, such as Movement Ecology as in G.H.  Weiss \cite{Weiss}, J. Masoliver and G.H. Weiss \cite{Masoliver}, K.P. Hadeler \cite{Hadeler1996,Hadeler1999}, T. Hillen \cite{Hillen2010}, and V. Méndez, D. Campos, and F. Bartumeus \cite{Bartumeus}, among others. The system can be written as
\begin{equation}\label{main0}
\begin{cases}
u_t+\gamma u_x=\mu(-u+v)\hspace{1cm} x\in(0,L)\\
v_t-\gamma v_x=\mu(u-v)\hspace{1cm} x\in(0,L)\\
u(0,t)=v(L,t)=0.
\end{cases}
\end{equation}
In this case, the variables $u(x,t)$ and $v(x,t)$ represent the probability density of a single individual to be at the point $x$ at time $t$ moving respectively right or left with velocity $\gamma>0$. Individuals can change direction depending on the difference $u(x,t)-v(x,t)$ with a rate $\mu>0$ (see \cite{Kac,Hadeler1996,Hadeler1999}). The individual leaves the system when it arrives at some of the ends. In some biological models, it leaves the system because it has found there the food it was looking for.

An interesting nonlinear modification of this model (see \cite{Hadeler1996,Hadeler1999}) is
\begin{equation*}
\begin{cases}
u_t+\gamma u_x=\mu(-u+v)+F(u+v)\hspace{1cm} x\in(0,L)\\
v_t-\gamma v_x=\mu(u-v)+F(u+v)\hspace{1cm} x\in(0,L),
\end{cases}
\end{equation*}
that admits homogeneous equilibrium solutions $u_0(x,t)\equiv v_0(x,t)\equiv c$ when $F(2c)=0$. In this case, when adding the boundary conditions $u(0,t)=v(L,t)=c$, the linearized system around this homogeneous equilibrium is of particular interest, and if $F'(2c)=a$ it can be written as
\begin{equation}\label{mainlin}
\begin{cases}
u_t+\gamma u_x=-(\mu+a)u+(\mu+a)v+2au\hspace{1cm} x\in(0,L)\\
v_t-\gamma v_x=(\mu+a)u-(\mu+a)v+2av \hspace{1cm} x\in(0,L)
\end{cases}
\end{equation}
with boundary conditions $u(0,t)=v(L,t)=0$.
The previous system has the same form as \eqref{main0}, as long as $\mu+a>0$ (except that the eigenvalues have been translated by the addition of the constant quantity $2a$). Therefore, the analysis of eigenvalues and eigenfunctions of \eqref{main0} that will be done in the present paper will remain the same, except for the translation. And this is very relevant for the stability of such equilibria.

Another interest in studying system \eqref{main0} as the linearization around a homogeneous positive state is that the non-positive solutions or eigenfunctions of \eqref{mainlin} will still have physical meaning, being understood as small perturbations of the homogeneous state.

In the new variables $t'=\mu t$ and $x'=x/L-1/2$ and defining the non-dimensional parameter $S=\gamma/(\mu L)>0$ system \eqref{main0} can be written in the following non-dimensional form:
\begin{equation}\label{main}
\begin{cases}
u_t+S u_x=-u+v\hspace{1cm} x\in(-1/2,1/2)\\
v_t-S v_x=u-v\hspace{1cm} x\in(-1/2,1/2)\\
u(-1/2,t)=v(1/2,t)=0,
\end{cases}
\end{equation}
where we have written $t$ and $x$ instead of $t'$ and $x'$.

It is usual to consider the previous system as a model in Movement Ecology as we have said above, but other applications are possible. For example, to model the phenomenon of durotaxis, which consists of the movement of certain biological cells towards the stiffer parts of substrate tissue (see C.R. Doering et al.\cite{Doering}), or to model the spatial distribution of growing cytoskeletal elements (see A. Büttenschön and L. Edelstein-Keshet in \cite{Buttenschon}). Also, this system was a starting point for other models of chemotaxis and self-organized biological aggregations (F. Lutscher \cite{Lutscher2002}, R. Eftimie \cite{Eftimie2012}).

But it is interesting to notice that the same system also appears in such a different field as Chemical Engineering or Biotechnology as a model of countercurrent parallel-flow heat exchangers (J.H. Chen and L. Malinowski \cite{ChenMalinowski}), and in chromatography as a linear model of countercurrent operation for species separation (see V.T.M. Silva et al \cite{Silva2004}, F. Wei et al. \cite{Wei}, V. Grosfils et al. \cite{Grosfils2007}, P. Satzer et al. \cite{Satzer2014}, Y. Guo et al. \cite{Guo2023}).

More concretely, in J. Menacho and J. Solà-Morales \cite{Menacho} the state of a single chromatography column of a true-moving bed adsorber is represented by the solute concentrations $c(x,t)$ in the liquid phase, and $q(x,t)$ in the solid phase in the interval $-L/2<x<L/2$. Then it is considered as boundary conditions that this single column is being washed out by the injection of pure solvent through the entrance port $x=-L/2$ in the liquid phase, and the injection of pure adsorbent at $x=L/2$ in the solid phase. Then, these concentrations evolve according to the following system of equations:
\begin{equation*}
\begin{cases}
c_t+v_\ell c_x=-FK(Hc-q) \hspace{1cm} x\in(-L/2,L/2)\\
q_t-v_s q_x=K(Hc-q)\hspace{1cm} x\in(-L/2,L/2)\\
c(-L/2,t)=0,\  q(L/2,t)=0,
\end{cases}
\end{equation*}
where all the parameters are positive: $v_\ell,v_s$ are the liquid and solid phase velocities, respectively, $H$ is the phase equilibrium constant, $F$ is the volume ratio between the solid and liquid phases, and $K$ is the adsorption kinetic constant.
We consider the simplified version of this system, which can be called {\sl reversible}, where the roles of $c(x,t)$ and $q(x,t)$ can be interchanged. That is, if $(c(x,t),q(x,t))=(u(x,t),v(x,t))$ is a solution, then $(c(x,t),q(x,t))=(v(-x,t),u(-x,t))$ is also a solution. This happens only when $v_\ell=v_s=v$ and $F=H=1$. The problem would then become
\begin{equation*}
\begin{cases}
c_t+vc_x=-K(c-q) \hspace{1cm} x\in(-L/2,L/2)\\
q_t-vq_x=K(c-q) \hspace{1cm} x\in(-L/2,L/2)\\
c(-L/2,t)=q(L/2,t)=0.
\end{cases}
\end{equation*}
By using the non-dimensional variables $x'=x/L$, $t'=Kt$, we obtain system \eqref{main}
with $S=v/(KL)$. As we will see, our goal in this work is to see the exponential rate at which these solutions approach the unique equilibrium $c_0(x,t)\equiv q_0(x,t)\equiv 0$, that is, how long will it take to be completely cleaned off. Observe that if in these new variables, we obtain a decay like $e^{-\alpha t}$, then in the original variables it would be $e^{-\alpha Kt}$. In fact, the present work partially completes the previous one in \cite{Menacho}, in the sense that now we will be able to obtain the exact rate of decay of the solutions.

\begin{rem}
When we consider system \eqref{main} as a model in chromatography, the application interest is for small values of the parameter $S>0$. Indeed, taking $S>1$ would mean that the transfer velocity between the liquid and solid phases would be much smaller than the velocity at which these phases move, making the first one negligible in practice. In ecology applications, however, large values of $S>0$ would mean a very fast emigration (compared to the number of direction changes), and this could be of practical interest.
\end{rem}

System \eqref{main} is strongly related to the telegraph equation, also known as the weakly damped wave equation. The first relation, and the most usually present in the bibliography, is the following one. According to \cite{Weiss,Masoliver, Hadeler1996,Hadeler1999}, system \eqref{main} can be transformed to an equivalent system for the total population density $p=u+v$, and the population flow $q= (u - v)$:
\begin{equation*}
\begin{cases}
p_t+Sq_x=0 \\
q_t+S p_x=-2 q.
\end{cases}
\end{equation*}
Taking derivatives in the first equation with respect to $t$, and in the last equation with respect to $x$, combining both new equations, and observing that the variable $q$ can be eliminated using that $p_t=-Sq_x$, we arrive at the following equation, known as the telegraph equation, or also weakly damped wave equation,
\begin{equation}\label{eq:onesp}
 p_{tt} + 2p_t - S^2 p_{xx}=0.
\end{equation}
Concerning the boundary conditions, it can be seen that they are translated into the following dynamic boundary conditions
\begin{equation}\label{eq:cconesp}
p_t(-1/2,t)-Sp_x(-1/2,t)+2p(-1/2,t)=0, \ \textrm{and }  p_t(1/2,t)+Sp_x(1/2,t)+2p(1/2,t)=0
\end{equation}
(see formula (60) in \cite{Hadeler1996}). Observe that following a very similar process we can see that $q$ also satisfies the telegraph equation \eqref{eq:onesp}, now with boundary conditions
\begin{equation}\label{eq:cconesq}
q_t(-1/2,t)-Sq_x(-1/2,t)=0, \ \textrm{and } q_t(1/2,t)+Sq_x(1/2,t)=0.
\end{equation}
If one repeats the previous procedure for the dimensional form of the equation (that is \eqref{main0}), one gets $ \frac{1}{2\mu}\,p_{tt} + p_t -\frac{\gamma^2}{2\mu}\, p_{xx}=0$.  Letting the parameters $\mu$ and $\gamma$ go to infinity such that the quotient $\frac{\gamma^2}{2\mu}$ converges to a finite value $D > 0$, we obtain, in the formal limit, the parabolic diffusion equation $p_t = D p_{xx}$.

But there is a second relation between \eqref{main} and the weakly damped wave equation. Indeed, we can see that if $u$ and $v$ satisfy system \eqref{main} then both $u$ and $v$ also satisfy the scalar weakly damped wave equation
\begin{equation}\label{dw}
w_{tt}+ 2w_t-S^2w_{xx}=0.
\end{equation}
This is a well-known fact (maybe not so usual in the literature as the previous one), but it is worth to be discussed here.
Applying the differential operator $\partial_t-S\partial_x$ to the first equation in \eqref{main} one gets $u_{tt}-S^2u_{xx}=(v_t-Sv_x)-u_t+Su_x$, and using the second equation $u_{tt}-S^2u_{xx}=(u-v)-u_t+Su_x$. Using now the first equation again one has $u_{tt}-S^2u_{xx}=(-u_t-Su_x)-u_t+Su_x=-2u_t.$

Concerning the boundary conditions one has of course $u(-1/2,t)=0$, but the condition at $x=1/2$ becomes more complicated: the only reasonable thing seems to be to evaluate at $x=1/2$ the terms of the first equation in \eqref{main}, use that $v(1/2,t)=0$, and obtain a dynamic boundary condition. Therefore, the boundary conditions are:
\begin{equation}\label{bc}
u(-1/2,t)=0,\ \textrm{and} \ u_t(1/2,t)+Su_x(1/2,t)+u(1/2,t)=0.
\end{equation}

Something similar happens with the unknown $v(x,t)$, except that the boundary condition at $x=-1/2$ is slightly different. In this case:
\begin{equation}\label{bc2}
v_t(-1/2,t)-Sv_x(-1/2,t)+v(-1/2,t)=0,\  \textrm{and } v(1/2,t)=0.
\end{equation}

The reverse property is now very easy. Suppose that $u(x,t)$ satisfies the wave equation \eqref{dw} together with the boundary conditions \eqref{bc}. Then, define $v(x,t)=u_t(x,t)+Su_x(x,t)+u(x,t)$. By this definition, it is clear that $u$ and $v$ satisfy the first equation in \eqref{main}. It is also clear that $v(1/2,t)=0$ and also, $(\partial_t-S\partial_x)v=u_{tt}-S^2u_{xx}-(\partial_t-S\partial_x)u=-2u_t+u_t-Su_x=-u_t+Su_x=u-v$, the second equation in \eqref{main}. A similar process allows us to prove the bi-directionality between \eqref{main} and \eqref{dw} and \eqref{bc2} for $v$.

It has to be said that, although it can seem that both problems are uncoupled and, hence, one can obtain $u$ or $v$ independently, they are indeed \textsl{weakly} coupled through the necessary initial conditions that relate one variable to the other one. For instance, $u_t(0,x)= -Su_x(0,x)-u(0,x)+v(0,x)$.

The difference between these two relations between system \eqref{main} with the weakly damped wave equation mainly relies on two facts. First, in \eqref{eq:onesp}-\eqref{eq:cconesp} and \eqref{eq:onesp}-\eqref{eq:cconesq} (the approach used in \cite{Hadeler1996}, among others) one is dealing with $p=u+v$ and $q=u-v$, not with $u,v$ separately. In this case, the equivalence is not completely straightforward (see the comment on the one-parameter family of solutions mentioned on p.138 of \cite{Hadeler1996}), and it may depend on the boundary conditions. But the relation and the physical interpretation between the solutions of system \eqref{main} and those of the wave equations \eqref{dw}-\eqref{bc} and \eqref{dw}-\eqref{bc2} is direct in our case, as we have just seen. Second, observe that the dynamic boundary conditions \eqref{eq:cconesp} for $p$ or \eqref{bc} and \eqref{bc2} for $u,v$ are not exactly the same.

The equivalence between system \eqref{main} and the weakly damped wave equation (in both approaches) seems a very relevant property, but we recognize that we have not used it in the present paper. The weakly damped wave equation with Dirichlet boundary conditions has been widely studied. A recent work on damped wave equations with dynamic boundary conditions but with strong damping instead of weak one is N. Fourrier and I. Lasiecka \cite{FourLas}, for instance. However, to our knowledge, the weakly damped wave equation with dynamic boundary conditions \eqref{bc} or \eqref{bc2} (or even \eqref{eq:cconesp} or \eqref{eq:cconesq}) are not very considered in the literature, and they are \textsl{equivalent} to system \eqref{main} in the sense described above. From our point of view, this also adds more interest to the present work.

As we said, dissipative wave equations are a very active research field nowadays. We will see in Section \ref{sec:vapsfups} that system \eqref{main} behaves similarly: the solutions oscillate in space and time, move and, finally, dissipate. Moreover, we will see that we have Sturm-Liouville-type space oscillations, in the sense that the $n$-th eigenfunction makes $n$ complete half-turns in the complex plane (see Proposition \ref{prop:osc}). In particular, as the long-time behaviour is given by the first eigenvalue and the corresponding eigenfunction, the solutions oscillate less and less as time increases. This is a behaviour that was first observed in the heat equation by Sturm: as time increases, the quick oscillatory and dissipative effect of higher eigenvalues (and eigenfunctions) tends to decrease, finally ending with the asymptotic decay rate and the non-oscillatory asymptotic profile given by the dominant eigenvalue and the corresponding eigenfunction.

Our goal will be to study the long-time behaviour of the solutions of system \eqref{main}.
For that purpose, it will be convenient to write it in an abstract form. So, we consider the function space $X=L^2(-1/2,1/2)\times L^2(-1/2,1/2)$ with the usual inner product, and write \eqref{main} as the following first order evolution equation:
\begin{equation}\label{opA}
\dfrac{d}{dt}U(t) = A U(t) ,\ \ t\in[0,+\infty)\vspace{0.2cm}
\end{equation}
where $U(t)=(u,v)^T$, and $A:\cD(A)\subset X\longrightarrow X$ is the linear operator
\begin{equation*}
A =
\begin{pmatrix}
-S\frac{d}{dx}-Id&Id\\Id&S\frac{d}{dx}-Id
\end{pmatrix}
\end{equation*}
with domain
$$\cD(A)=H^1_l(-1/2,1/2) \times H^1_r(-1/2,1/2).$$
Here, $H^1_l(-1/2,1/2)$ and $H^1_r(-1/2,1/2)$ are the subspaces of the Sobolev space $H^1(-1/2,1/2)$ of the functions that vanish, respectively, at the left and the right end of the interval.

One can see that problem \eqref{opA} is well-posed. Actually, this is proved in Theorem 3.4 of T. Hillen \cite{Hillen2010}, which we include here for the sake of completeness of the present work.
\begin{thm}\label{thm:reg}[Theorem 3.4 of \cite{Hillen2010}]
The operator $A$ is the generator of a $C_0$-semigroup on $X$, $\{ T(t)=e^{At}, \ t\geq 0\}$. In particular, for any $U_0\in \cD(A)$, there exists a unique solution $U=T(t)U_0\in C^1([0,\infty);X))\cap C^0([0,\infty);\cD(A))$ satisfying \eqref{opA}. If $U_0\in X$, then $U=T(t)U_0$ is the unique solution in $\in C^0([0,\infty);X)$ satisfying \eqref{opA} in the mild sense.
\end{thm}

As we said, Theorem \ref{thm:reg}  is a result given in \cite{Hillen2010}. Actually, he proves that $A$ is the generator of a $C_0$- semigroup on $L^p(-1/2,1/2) \times L^p(-1/2,1/2)$, for $1\leq p <\infty$. For the purpose of the current work, $p=2$ suffices.

A $C_0$-semigroup $T(t)$ on a Banach space $X$ is called eventually compact if there exists $t_0>0$ such that $T(t)$ is compact for $t>t_0$ (see K.J. Engel and R. Nagel \cite{EngelNagel}, for example).  The next proposition shows that our semigroup becomes compact for $t>1/S$. This is consistent with the ideas in  \cite{Hadeler1999} Sect. 4 and in \cite{Hillen2010}, that the irregularities of the initial conditions are washed out to the boundary after a time $t=1/S$.

\begin{thm}\label{thm:compact}
The operator $(A, D(A))$ given in \eqref{opA} is the infinitesimal generator of an eventually compact semigroup $T(t)$ with $t_0=1/S$.
\end{thm}
The proof of this result will be given in Section \ref{sec:decay}, and will be a key point in the proof of the decay rate of the solutions of our problem (given in Theorem \ref{thm:main} below).

The aim of this work is to prove that the solutions of the problem \eqref{main} decay exponentially in the norm of $X$, and this decay is optimal; that
this optimal decay rate is given by the dominant eigenvalue of $A$, which can be calculated as a particular solution of a system of transcendental equations; and that the solutions tend to 0 with a certain asymptotic profile. These results are given in the following theorem, which summarizes the main results of the present work.

\begin{thm}\label{thm:main}
For the correlated random walk system \eqref{main} (with abstract form given in \eqref{opA}), the following assertions are true:
\begin{enumerate}
\item For each $S>0$ there exists a real number $\lambda_0(S) < 0$ and a real number $M\geq 1$ such that
$$\| e^{At}\| \leq M e^{\lambda_0(S) t} , \textrm{ for all } t\geq 0$$
where the norm is the usual one in the space $X$ where the semigroup $e^{At}$ is defined.

\item The previous decay is optimal. That is, there exists a solution of \eqref{main} (the eigenfunction whose eigenvalue is $\lambda_0(S)$) such that decays exactly at the previous rate.

\item For all $S>0$, this number  $\lambda_0(S)<0$ is the dominant eigenvalue of $A$. That is, all the other eigenvalues $\lambda$ of $A$ satisfy the inequality and gap condition
        $${\Re}(\lambda(S))<\lambda_0(S)-\varepsilon$$
      for some $\varepsilon=\varepsilon(S)>0$. Also, $\lambda_0(S)$ has geometric and algebraic multiplicity equal to 1.

\item For an open and dense set of initial conditions, the solutions have an asymptotic profile when $t\to\infty$ that is a multiple of the eigenfunction corresponding to the eigenvalue $\lambda_0(S)$, $(u_0(x),v_0(x))$, in the sense that they tend to multiples of this function at a faster velocity than the one at which they tend to $0$. This means that they turn into the eigenfunction before becoming $0$.
\end{enumerate}
\end{thm}

The \textsl{explicit} description of the dominant eigenvalue $\lambda_0(S)$ and its corresponding eigenfunction $(u_0(x),v_0(x))$ is given in point 5 of Proposition \ref{prop:descvaps} and Remark \ref{rem:domfup}, respectively. It has to be said that some of the results of point 5 of Proposition \ref{prop:descvaps} can be found in \cite{Hadeler1996,Hadeler1999}, although there the dominant eigenvalue $\lambda_0$ is given as a solution of a characteristic equation of an equivalent problem, and also given in a different form (although equivalent) to ours. In this sense, we can say that we resume and complete the spectral analysis started by K.P. Hadeler in these works, as we also describe the rest of the eigenvalues, and eigenfunctions of this problem.

Comparing our approach with the one in the works of K.P. Hadeler, we would say that in our case we formulate the problem in a way that highlights the symmetries of the problem. To begin with, system \eqref{main} is posed in $(-1/2,1/2)$ instead of $(0,L)$ (as in \cite{Hadeler1996}). Also, we will see that the fact of working with $u,v$ (instead of $p=u+v, q=u-v$ as in \cite{Hadeler1996}), and the way in which we write and compute the eigenvalues and eigenfunctions of \eqref{main} makes explicit the symmetric and antisymmetric character of the eigenfunctions, that remains hidden in the $p,q$ approach (see Remark \ref{rem:sym}). Also, the way in which we ultimately compute the eigenvalues, with an auxiliary parameter $\nu$, substantially simplifies the characterization of the spectrum  (see Propositions \ref{prop:vapsfups} and \ref{prop:descvaps})

\begin{rem}
The results in Theorem \ref{thm:main} above imply that $S>0$ is a dissipative parameter in the following sense: for $S$ large, the solutions tend to $0$ faster than for smaller $S$ (see also, for instance, Figure \ref{fig:ReLambda}).
\end{rem}

\begin{rem}
A complete description of $\lambda_0(S)$ and $(u_0(x),v_0(x))$, together with the rest of the spectrum $\lambda_{n,j}(S)$ and eigenfunctions $\left( u_{n,j}(x), v_{n,j}(x)\right)$, $n\geq 1,\, j=1,2$, is given in Lemma \ref{lem:spec}, and Propositions \ref{prop:vapsfups} and \ref{prop:descvaps}, in Section \ref{sec:vapsfups} below. Indeed, we will see that all the eigenfunctions and the eigenvalues (in particular the dominant ones, $\lambda_0(S)$ and $(u_0(x),v_0(x)$) are solution of an explicit and relatively simple transcendental system of equations.
\end{rem}

Similar spectral analysis techniques were used by part of the authors of the present work in the recent paper \cite{PellSM2019}, where also a comprehensive and accurate spectrum description was necessary to obtain the optimal exponential decay rate of the solutions of a Moore-Gibson-Thompson equation.

In the present work, this complete description of the eigenfunctions and eigenvalues is the first step in the proof of Theorem \ref{thm:main}, and is given in Section \ref{sec:vapsfups}, together with the Sturm-Liouville type space oscillations result. Once this description is done, we start with the proof of part 3 of Theorem \ref{thm:main}, that is, by proving that $\lambda_0(S)$ is, indeed, the dominant eigenvalue (see Section \ref{sec:domvap}). After this, we proceed with the proofs of parts 1, 2, and 4 of Theorem \ref{thm:main}, which are done in Section \ref{sec:decay}. In this section we also give the proof of the compactness result stated in Theorem \ref{thm:compact}, which turns out to be crucial to prove the optimal asymptotic decay rate of the solutions of the problem.

\section{Eigenfunctions and eigenvalues}\label{sec:vapsfups}

In this section, we give a complete description of the eigenvalues and eigenfunctions of the system whose abstract form is defined in \eqref{opA}.

\begin{lemma}\label{lem:spec}
The spectrum of the operator $A$ defined in \eqref{opA} consists only of eigenvalues (which will be described in Proposition  \ref{prop:descvaps}).
\end{lemma}

\begin{pf}
We observe that $A$ has compact resolvent, because of the compact embedding $D(A)=H^1_\ell(-1/2,1/2)\times H^1_r(-1/2,1/2)\subset L^2(-1/2,1/2)\times L^2(-1/2,1/2)=X$. Therefore, $A$ has only point spectrum.
\end{pf}

\begin{prop}\label{prop:vapsfups}{(Eigenfunctions and eigenvalues)}

Given $S>0$, we consider the following eigenvalue system
\begin{equation}\label{mainev}
\begin{cases}
-Su'(x)-u(x)+v(x)=\lambda u(x)\hspace{.5cm} -1/2<x<1/2,\\
Sv'(x)+u(x)-v(x)=\lambda v(x)\hspace{.5cm} -1/2<x<1/2,\\
u(-1/2)=0,\\
v(1/2)=0.
\end{cases}
\end{equation}
The solutions of this system give the eigenvalues $\lambda\in\mathbb{C}$ and the eigenfunctions $(u(x),v(x))^T$ of the operator $A$ defined in \eqref{opA}. Then, the following assertions hold.
\begin{enumerate}
\item System \eqref{mainev} admits nontrivial solutions in the following three cases, and only in these cases:
\begin{enumerate}
\item there exists a $\nu\in\C\setminus\{0\}$ such that
\begin{equation}\label{ces}
\begin{cases}
\sin(\nu)=S\nu,\\
\lambda=-1-\cos(\nu),
\end{cases}
\end{equation}
(and then $\lambda\ne-2$).
\item $\lambda=-2$, only when $S=1$
\item there exists a $\nu\in\C\setminus\{0\}$ such that
\begin{equation}\label{cea}
\begin{cases}
\sin(\nu)=-S\nu,\\
\lambda=-1+\cos(\nu)
\end{cases}
\end{equation}
(and then $\lambda\ne0$).
\end{enumerate}
In the cases (a) and (c) we have $\nu=\sqrt{-\lambda^2-2\lambda}\ /S$, but both signs of the square root are admissible, as they will give the same eigenfunction (see Remark \ref{rem:sqrt}).

\item The previous eigenvalues are geometrically simple, that is the space of solutions has dimension one. The corresponding eigenfunction in each previous case is:
\begin{enumerate}
\item
\begin{equation}\label{eigenforms_sym}
\begin{pmatrix} u(x)\\v(x)\end{pmatrix}=
\begin{pmatrix} \sin(\nu(1/2+x))\\\sin(\nu(1/2-x))\end{pmatrix}
\end{equation}

\item
\begin{equation}\label{eigenforms_-2}
\begin{pmatrix} u(x)\\v(x)\end{pmatrix}=
\begin{pmatrix} 1+2x\\1-2x\end{pmatrix}
\end{equation}
\item
\begin{equation}\label{eigenforms_asym}
\begin{pmatrix} u(x)\\v(x)\end{pmatrix}=
\begin{pmatrix} \sin(\nu(1/2+x))\\-\sin(\nu(1/2-x))\end{pmatrix}.
\end{equation}

\end{enumerate}
\end{enumerate}

\end{prop}

\begin{figure}[htpb]
    \centering
   \includegraphics[width=0.3\textwidth]{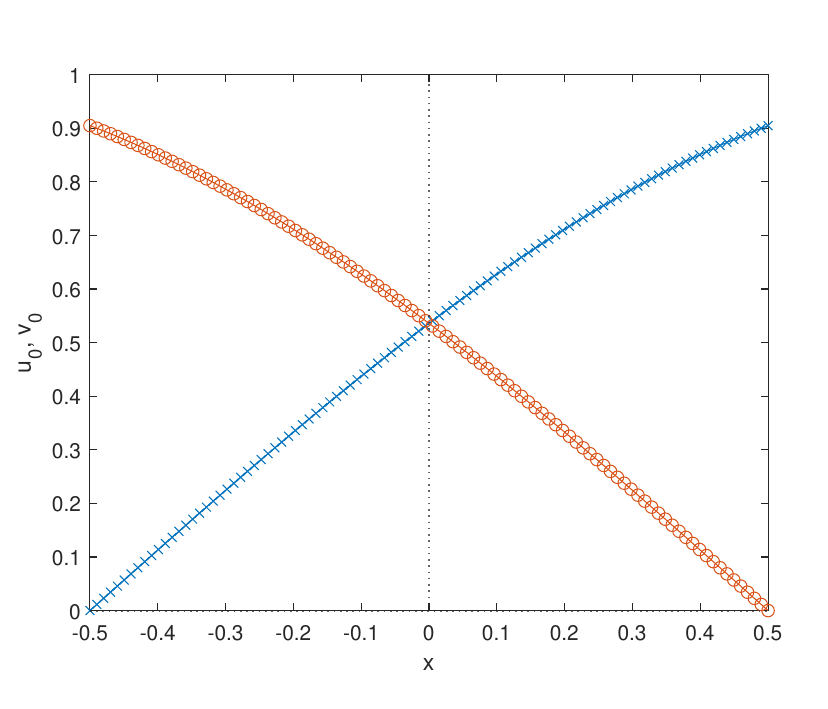}\hspace{0.5cm}
   \includegraphics[width=0.3\textwidth]{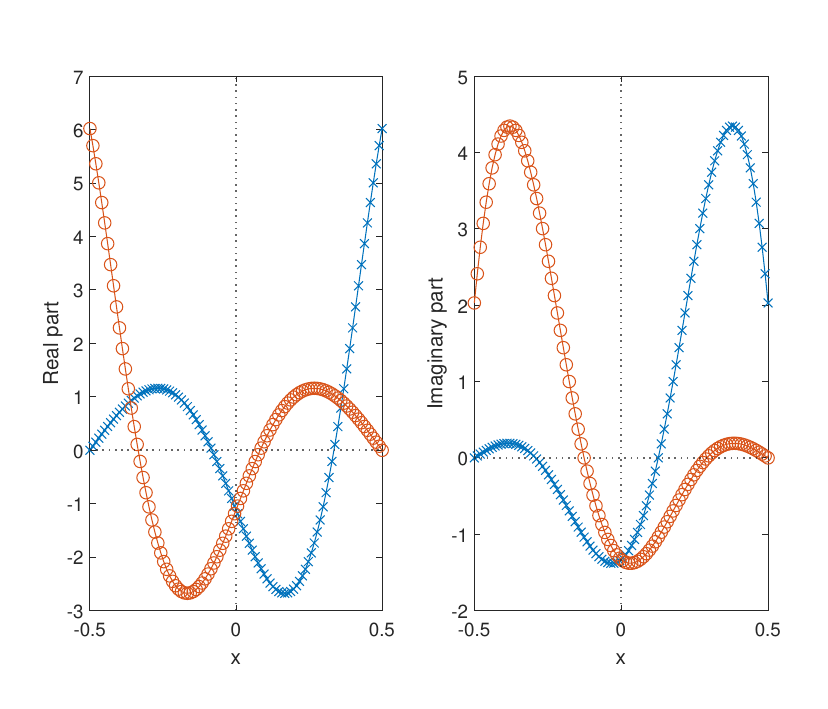}\hspace{0.5cm}
   \includegraphics[width=0.3\textwidth]{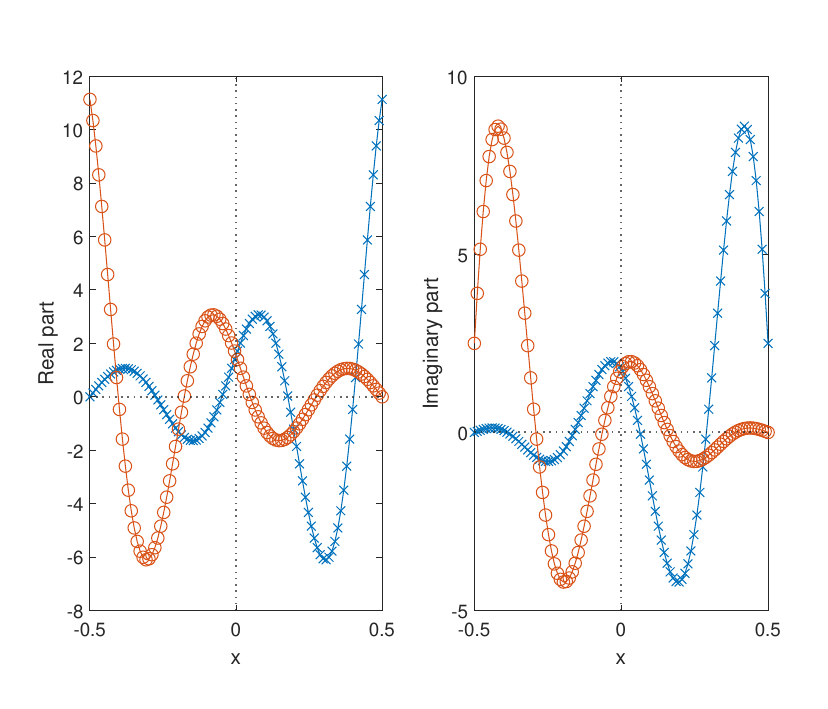}
   \includegraphics[width=0.3\textwidth]{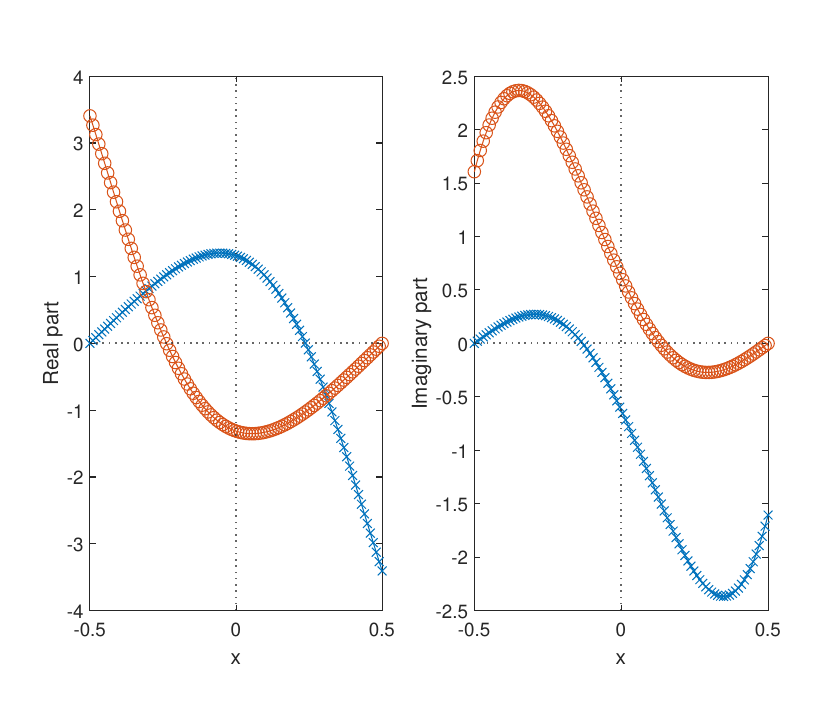}\hspace{0.7cm}
   \includegraphics[width=0.3\textwidth]{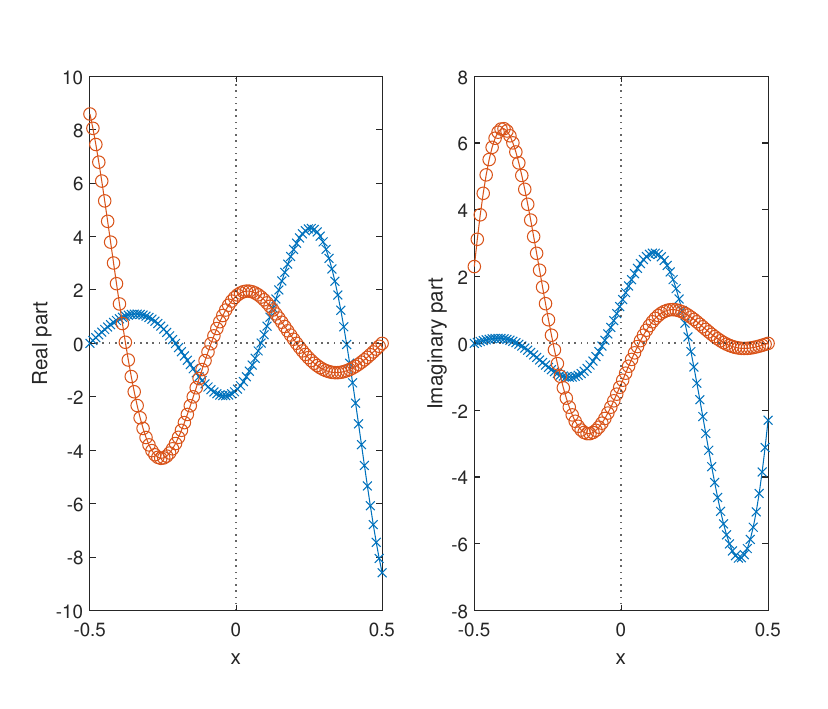}
   \caption{Some eigenfunctions for $S=0.8$, real and imaginary parts respect to $x\in[-1/2,1/2]$.  From left to right, and up to down, plots from $u_0,v_0$, and real and imaginary part of $u_{n,1},v_{n,1}$ for $n=2,4$ (symmetric case, top), and $n=1,3$ (antisymmetric case, bottom). In all cases, the blue-cross line corresponds to $u_{n,1}$, and the red-circle one to $v_{n,1}$. The plots for the corresponding $(u_{n,2},v_{n,2})$ are similar, so we do not include them. In these graphs, we can also see the oscillatory behaviour of the $n$-the eigenfunction, described in Proposition \ref{prop:osc} below. }.
    \label{fig:FUPS}
\end{figure}

\begin{rem}\label{rem:sym}
In cases (a) and (b) the corresponding eigenfunctions can be called of a {\sl symmetric} type because they are of the form $(u(x),u(-x))$. In case (c), they can be called of {\sl antisymmetric} type since it has the form $(u(x),-u(-x))$. See Figure \ref{fig:FUPS}.
\end{rem}

\begin{rem}\label{rem:sqrt}{(Notation Remark)}
For a complex number $z$, to avoid the ambiguity of $\sqrt{z}$, we will use from now on the notation $\sqr(z)$ to denote the unique square root of $z$ that satisfies $\sqr(0)=0$ and $-\pi/2 <\arg(\sqr(z))\le \pi/2$.
\end{rem}

\begin{pf}
The system of equations for the eigenfunctions of eigenvalue $\lambda$ of the operator $A$ defined in \eqref{opA} is
\begin{equation}\label{lambdaeq}
\dfrac{d}{dx}
\begin{pmatrix} u\\ v\end{pmatrix}=
\begin{pmatrix}
(-1-\lambda)/S&1/S\\-1/S&(1+\lambda)/S
\end{pmatrix}\begin{pmatrix} u\\ v\end{pmatrix},
\end{equation}
where, in general, $\lambda$ will be a non-real complex number.

For $\lambda\ne 0,-2$ let us define the auxiliary parameter $\nu=\nu(\lambda,S)=\sqr(-\lambda^2-2\lambda)/S$ (see Remark \ref{rem:sqrt}). Observe that the excluded values $\lambda= 0,-2$ are precisely the values such that $-\lambda^2-2\lambda=0$.

For $\lambda\ne 0,-2$, it is easy to see that the general solution of \eqref{lambdaeq} is
\begin{equation}\label{eigen}
\begin{pmatrix}u\\v\end{pmatrix}=
C_1e^{i\nu x}
\begin{pmatrix} 1\\ 1+\lambda+iS\nu\end{pmatrix}
+
C_2e^{-i\nu x}\begin{pmatrix} 1+\lambda+iS\nu\\1\end{pmatrix}.
\end{equation}
Imposing the boundary conditions $u(-1/2)=v(1/2)=0$ we get the system of equations
\begin{equation}\label{Cs}
\begin{cases}
C_1e^{-i\nu/2}+C_2e^{i\nu/2}(1+\lambda+iS\nu)=0,\\
$$C_1e^{i\nu/2}(1+\lambda+iS\nu)+C_2e^{-i\nu/2}=0,
\end{cases}
\end{equation}
that in order to have nontrivial solutions requires
$e^{i\nu}(1+\lambda+iS\nu)^2=e^{-i\nu}, $ or, equivalently,
\begin{equation}\label{pm}
(1+\lambda+iS\nu)=\pm e^{-i\nu}.
\end{equation}

By taking the  $-$ sign in \eqref{pm}, the two equations in system \eqref{Cs} become the single equation $C_1e^{-i\nu/2}-C_2e^{-i\nu/2}=0$, that implies $C_1=C_2$. In this case, the eigenfunction given by \eqref{eigen} is
\begin{equation*}
\begin{pmatrix}u\\v\end{pmatrix}=
C_1e^{i\nu x}
\begin{pmatrix} 1\\ -e^{-i\nu}\end{pmatrix}
+C_1e^{-i\nu x}\begin{pmatrix} -e^{-i\nu}\\1\end{pmatrix}
=2iC_1e^{-i\nu/2}\begin{pmatrix}
\sin(\nu(1/2+x))\\ \sin(\nu(1/2-x))\end{pmatrix},
\end{equation*}
and no other non-trivial solutions of \eqref{Cs} exist with $(1+\lambda+iS\nu)=- e^{-i\nu}$.

Now, with the $+$ sign in \eqref{pm} the two equations in system \eqref{Cs} become the single equation $C_1e^{-i\nu/2}+C_2e^{-i\nu/2}=0$, that implies $C_1=-C_2$. In this case, the eigenfunction given by \eqref{eigen} is
\begin{equation*}
\begin{pmatrix}u\\v\end{pmatrix}=
C_1e^{i\nu x}
\begin{pmatrix} 1\\ e^{-i\nu}\end{pmatrix}
-C_1e^{-i\nu x}\begin{pmatrix} e^{-i\nu}\\1\end{pmatrix}
=2iC_1e^{-i\nu/2}\begin{pmatrix}
\sin(\nu(1/2+x))\\ -\sin(\nu(1/2-x))\end{pmatrix},
\end{equation*}
and no other solutions of \eqref{Cs} exist with $(1+\lambda+iS\nu)= e^{-i\nu}$.

Summarizing, if $\lambda\ne 0, -2$ is an eigenvalue of \eqref{opA}, then there are two alternative possibilities for the eigenfunctions, and no more than these two: \eqref{eigenforms_sym} or \eqref{eigenforms_asym}, that is
\begin{equation*}
\begin{pmatrix} u\\v\end{pmatrix}=
\begin{pmatrix}
\sin(\nu(1/2+x))\\
\sin(\nu(1/2-x))\end{pmatrix}
\hbox{  or  else }
\begin{pmatrix} u\\v\end{pmatrix}=
\begin{pmatrix}
\sin(\nu(1/2+x))\\ -\sin(\nu(1/2-x))\end{pmatrix},
\end{equation*}
where $\nu=\nu(\lambda,S)=\sqr(-\lambda^2-2\lambda)/S$. The first one corresponds to the $-$ sign in \eqref{pm}, and the second one corresponds to the $+$ sign.

Let us consider now the exceptional cases $\lambda=0,-2$. If $\lambda=0$, the general solution of \eqref{lambdaeq} is
\begin{equation*}
\begin{pmatrix}u\\v\end{pmatrix}=
C_1
\begin{pmatrix} 1\\ 1\end{pmatrix}
+
C_2\begin{pmatrix} x\\S+x\end{pmatrix},
\end{equation*}
and no possibility of adjusting to $u(-1/2)=v(1/2)=0$ exists except when $C_1=C_2=0$. The only possibility would be when $S=-1$, which is out of our scope. Therefore, $\lambda=0$ will never be an eigenvalue for $S> 0$.

The case $\lambda=-2$ is different. In this case, the general solution of  \eqref{lambdaeq} is
\begin{equation}\label{-2}
\begin{pmatrix}u\\v\end{pmatrix}=
C_1
\begin{pmatrix} 1\\ -1\end{pmatrix}
+
C_2\begin{pmatrix} x\\S-x\end{pmatrix}.
\end{equation}
It is impossible to adjust the boundary conditions for a nontrivial solution when $S>0$ except when $S=1$. In this case, there is no such obstruction, and the eigenfunction corresponding to $\lambda=-2,\,S=1$ would be \eqref{eigenforms_-2}, that is
\begin{equation*}
\begin{pmatrix}u\\v\end{pmatrix}=
\begin{pmatrix} 1+2x\\ 1-2x\end{pmatrix}.
\end{equation*}

At this moment we can prove the geometric simplicity of the eigenvalues. The expression \eqref{pm} can be written more explicitly as $1+\lambda+iS\nu(\lambda,S)=\pm e^{-i\nu(\lambda,S)}$. Since we have taken care of defining $\nu(\lambda,S)$ in an unambiguous way, and as $ e^{-\nu(\lambda,S)}$ is never zero, the previous equality can be satisfied with the $+$ sign or with the $-$ sign, but not with the two signs at the same time for the same values of $\lambda$ and $S$. This means that only one of the two forms \eqref{eigenforms_sym} or \eqref{eigenforms_asym} will be possible for given $\lambda$ and $S$. And they are linearly independent, since $\nu(\lambda,S)\ne0$, once we have excluded the cases $\lambda=0,-2$. In the case $\lambda=-2,\, S=1$ it is easy to see that \eqref{-2} admits only the nontrivial solution $C_2=2C_1$. Hence, we have proved the geometric simplicity of the eigenvalues. We will discuss the algebraic multiplicity of the dominant eigenvalue in Section \ref{sec:domvap}.

Let us now analyse $\lambda$ the eigenvalues of \eqref{opA}. Our goal is to obtain a characteristic equation for them. We start with the symmetric case, that is \eqref{eigenforms_sym}. We substitute $u(x)=\sin(\nu(1/2+x))$ and $v(x)=\sin(\nu(1/2-x))$ into the first row of the system $A(u,v)^T=\lambda(u,v)^T$, where $A$ is given in \eqref{opA}, and we get
\begin{equation}\label{sym1}
-S\nu\cos(\nu(1/2+x))-\sin(\nu(1/2+x))+\sin(\nu(1/2-x))=\lambda \sin(\nu(1/2+x)).
\end{equation}
Now, writing $\sin(\nu(1/2-x))$ as $-\sin(\nu(1/2+x-1))$, and this last expression as $$-\sin(\nu(1/2+x))\cos(\nu)+\cos(\nu(1/2+x))\sin(\nu)$$ we get
$$-S\nu\cos(\nu(1/2+x))-(1+\lambda)\sin(\nu(1/2+x))-\sin(\nu(1/2+x))\cos(\nu)+\cos(\nu(1/2+x))\sin(\nu)=0,$$
or, rearranging terms,
$$(-S\nu+\sin(\nu))\cos(\nu(1/2+x))-(1+\lambda+\cos(\nu))\sin(\nu(1/2+x))=0.$$
Due to the linear independence between the functions $\cos(\nu(1/2+x))$ and $\sin(\nu(1/2+x))$ when $\nu\ne 0$, we conclude that system \eqref{ces} holds, that is:
\begin{equation*}
\begin {cases}
\sin\nu=S\nu\\ \lambda=-1-\cos\nu.
\end{cases}
\end{equation*}

Therefore, we have proved that if $\lambda\ne -2$ is an eigenvalue with a symmetric eigenfunction of the form \eqref{eigenforms_asym} then \eqref{ces} must be satisfied with $\nu=\nu(\lambda,S)=\sqr(-\lambda^2-2\lambda)/S$,  $\nu\neq 0$.

We also observe that $\nu=0,\,\lambda=-2$ is a solution of system \eqref{ces} for all values of $S$, but we know that it is a spurious solution unless $S=1$ (in the sense that it does not have a non-zero eigenfunction).

From the proof, we can observe that the reciprocal is also true. That is,
if $\lambda\in\C\setminus\{-2\}$ and there exists a $\nu\in\C\setminus\{0\}$ such that \eqref{ces} holds, then
then $\lambda$ will be an eigenvalue with a symmetric eigenfunction of the form given \eqref{eigenforms_sym}.

To continue the analysis of the eigenvalues $\lambda$ of \eqref{opA}, we proceed now, in a similar way, with the antisymmetric case, that is \eqref{eigenforms_asym}. We substitute $u(x)=\sin(\nu(1/2+x))$ and $v(x)=-\sin(\nu(1/2-x))$ into the first row of $A(u,v)^T=\lambda(u,v)^T$ and we get
$$-S\nu\cos(\nu(1/2+x))-\sin(\nu(1/2+x))-\sin(\nu(1/2-x))=\lambda \sin(\nu(1/2+x)).$$
Now we write again $\sin(\nu(1/2-x))$ as $-\sin(\nu(1/2+x-1))$ and this last expression as $$-\sin(\nu(1/2+x))\cos(\nu)+\cos(\nu(1/2+x))\sin(\nu).$$
And substituting we get
$$-S\nu\cos(\nu(1/2+x))-(1+\lambda)\sin(\nu(1/2+x))+\sin(\nu(1/2+x))\cos(\nu)-\cos(\nu(1/2+x))\sin(\nu)=0,$$
or, rearranging terms,
$$(-S\nu-\sin(\nu))\cos(\nu(1/2+x))-(1+\lambda-\cos(\nu))\sin(\nu(1/2+x))=0.$$
Due again to the linear independence between the functions $\cos(\nu(1/2+x))$ and $\sin(\nu(1/2+x))$ when $\nu\ne 0$, we obtain system \eqref{cea}, that is:
\begin{equation*}
\begin {cases}
\sin\nu=-S\nu\\ \lambda=-1+\cos\nu.
\end{cases}
\end{equation*}
We have proved that if $\lambda$ is an eigenvalue with an antisymmetric eigenfunction of the form \eqref{eigenforms_asym}, then $\lambda\ne 0$  and \eqref{cea} is satisfied with $\nu=\nu(\lambda,S)=\sqr(-\lambda^2-2\lambda)/S$.

And we observe again the existence of a spurious solution of \eqref{cea}, $\nu=0,\ \lambda=0$, for all values of $S>0$.

From the proof, we can also observe that the reciprocal property is also true. That is, if $\lambda\in\mathbb{C}\setminus\{0\}$ and there exists a $\nu\in\C\setminus\{0\}$ such that equation \eqref{cea} is fulfilled, then
 $\lambda$ is an eigenvalue with an antisymmetric eigenfunction of the form \eqref{eigenforms_asym}.

With all this, we have proved all the parts of Proposition \ref{prop:vapsfups}.
\end{pf}

Let us now proceed to a detailed analysis of the solutions of the characteristic equations \eqref{ces} and \eqref{cea}.

\begin{prop}{(Description of the eigenvalues)}\label{prop:descvaps}\\
The eigenvalues $\lambda=\lambda(S)\in\mathbb{C}\setminus\{-2\}$ of \eqref{opA}, and the values $\nu=\nu(S)=\sqr(\lambda^2-2\lambda)/S\neq 0$, $S>0$,  given in Proposition \ref{prop:vapsfups} satisfy the following properties (see also Figures \ref{fig:nu} and \ref{fig:ReLambda} below).

\begin{enumerate}
\item For each $S>0$, the values of $\nu$, come in pairs (except the first one) and constitute an infinite sequence $\left\{ \nu_{0},\, \,\nu_{1,1},\,\nu_{1,2},\,\nu_{2,1},\, \nu_{2,2},\ldots  \nu_{n,1}, \,\nu_{n,2},\ldots\right\}\subset\mathbb{C}$, where $n$ even corresponds to solutions of equation \eqref{ces} and $n$
odd corresponds to solutions of \eqref{cea}.

\item There exists a decreasing sequence of values of $S$, say  $\{S_m\},m\geq 1$, with $0<S_m < 1$, satisfying that the equation $\sin(\nu)=S_m \nu$ ($m$ even) or the equation $\sin(\nu)=-S_m \nu$ ($m$ odd) has a double real root $\nu_m\in [m\pi,(m+1/2)\pi]$. Naturally, $S_m\to 0$ as $m\to\infty$. We call $S_{crit}:=S_1\simeq 0.2172$.

\item Regarding $\nu_{n,j}$, $n\geq 1$, $j=1,2$, we can say that (see Figures \ref{fig:nu} and \ref{fig:NUS1}):
\begin{enumerate}
\item For a fixed $0<S\leq S_{crit}$ with $S\in (S_{m+2},S_m]$, $m\geq 1$, the following assertions hold.
\begin{enumerate}
\item If $1\leq n < m$ (or if $n=m$ and $S\neq S_m$), we have $\nu_{n,j}=\nu_{n,j}(S)\in\mathbb{R},\ j=1,2$, with $\nu_{n,1}\neq \nu_{n,2}$. In this case,
$$\nu_{n,j}\in (n\pi,(n+1)\pi), \ j=1,2.$$

\item If $n=m$ and $S=S_m$, we have $\nu_{n,1}(S)=\nu_{n,2}(S)\in \mathbb{R}$ (double real root). In this case,
$$\nu_{n,j}\in (n\pi,(n+1/2)\pi), \ j=1,2.$$
\item If $n>m$, we have $\nu_{n,1}(S)=\overline{\nu_{n,2}(S)}\in \mathbb{C}\setminus\mathbb{R}$. In this case, $\nu_{n,1}\in Q_n$ and $\nu_{n,2}\in \overline{Q_n}$, where
\begin{equation}\label{eq:strip}
Q_n=\{z\in\mathbb{C}; \ n\pi < \Re(z) < (n+1/2)\pi,\ \Im(z)>0\}.
\end{equation}
Also, for $n\to\infty$ and $S>0$ fixed, we have
\begin{equation}\label{eq:nuasympt}
\nu_{n,j}(S)=(n+1/2)\pi-(-1)^j \log(2S(n+1/2)\pi)\,i+o(1),\ j=1,2.
\end{equation}
That is, to a leading order we have
$$\Re(\nu_{n,j}(S)) \simeq n\pi,\ \ \ \Im(\nu_{n,j}(S)) \simeq -(-1)^j\log(n),\ j=1,2.$$

\end{enumerate}
\item For a fixed $S>S_{crit}$, all $\nu_{n,j}(S)$, $n\geq 1$, $j=1,2$, the same conclusions as in point (iii) above hold.
\end{enumerate}

\item

Regarding $\lambda_{n,j}$, $n\geq 1$, $j=1,2$, we can say that (see also Figures \ref{fig:nu} and \ref{fig:ReLambda}):

\begin{enumerate}
\item In case (i), we have $\lambda_{n,1}(S)\neq \lambda_{n,2}(S)$, both in $\mathbb{R}$,  with
$$\lambda_{n,1}(S)=-1 -\sqr(1-S^2(\nu_{n,1}(S))^2),\ n\geq 1,$$
and
$$\lambda_{n,2}(S)=-1\pm \,\sqr(1-S^2(\nu_{n,2}(S))^2),\ n\geq 1.$$
The minus or plus (most common) sign of $\lambda_{n,2}$ depends on the parity of $n$ and the sign of $\cos(\nu_{n,2})$ (which can also be zero).
In particular, $\lambda_{n,1}(S)$ tend to $-2$ when $S\to 0$, and $\lambda_{n,2}(S)$ tend to $0$ when $S\to 0$.

\item In cases (ii) and (iii) (including (3b)), we have
\begin{equation*}
\lambda_{n,j}(S)=-1- \,\sqr(1-S^2(\nu_{n,j}(S))^2),\ n\geq 1,\ j=1,2.
\end{equation*}
(in case (ii) we have $\lambda_{n,1}(S) = \lambda_{n,2}(S)\in\mathbb{R}$).

\item In case (iii) (including (3b)) we have
\begin{equation}\label{lem:ineq}
\Re(\lambda_{n,j}(S))<-1-S.
\end{equation}

Also, for $S$ or $n$ large, we have
\begin{equation}\label{eq:lambdanu}
\lambda_{n,j}(S)= -1-\left((-1)^j S \nu_{n,j}(S)\, i\right)+o(1),\ j=1,2.
\end{equation}
In particular, when $n\to \infty$ and $S>0$ is fixed, we have
\begin{equation}\label{eq:lambdaasympt}
\lambda_{n,j}(S)=-1-S\log(2S(n+1/2)\pi)-(-1)^j S(n+1/2)\pi\,i+o(1),\ j=1,2
\end{equation}
which is,
$$\Re(\lambda_{n,j}(S)) \simeq -S\log(n),\ \ \ \Im(\lambda_{n,j}(S))\simeq -(-1)^j Sn\pi,\ j=1,2.$$

\end{enumerate}

\item The case $n=0$ is a bit different (see Figures \ref{fig:nu} and \ref{fig:ReLambda}).

First, we consider $\nu_0(S)$. In general, $\nu_0(S)$ is the solution of the transcendental equation $\sin(\nu)=S\nu$ whose real part is in $[0,\pi)$. Specifically,
\begin{enumerate}
\item if $0<S <1$, we have $\nu_0(S)\in\mathbb{R}$ and $\nu_0(S)\in(0,\pi)$. In the particular case $S=\frac{2}{\pi}$, we have $\nu_0(S)=\frac{\pi}{2}$;
\item if $S=1$, $\nu_0(S)=0$ ($\nu(S)=0$ is a spurious solution for all $S\neq 1$, since it does not have a non-zero eigenfunction);
\item if $S>1$, we have $\nu_0(S)\in i\mathbb{R}$, with $\Im (\nu_0(S)) >0$ (observe that we are not considering $\overline{\nu_0(S)}$, see Remark \ref{rem:sqrt}).
\end{enumerate}
Second, we consider $\lambda_0(S)$. Then, it is fulfilled that:
\begin{enumerate}
\item we have $\lambda_0(S)=-1-\cos(\nu_0(S))$. In particular, if $S=1$, we have $\lambda_0(S)=-2$, and if $S=\frac{2}{\pi}$, we have $\lambda_0(S)=-1$;
\item $\lambda_0(S)$ is real for all the values of $S$, and decreases in $S$;
\item $\lambda_0(S)\to 0$ as $S\to 0$, $\lambda_0(S)\to-\infty$ as $S\to\infty$.

\end{enumerate}
\end{enumerate}

\end{prop}

\begin{rem}\label{rem:domfup}
If we denote by $(u_0(x),v_0(x))$ the eigenfunction corresponding to the eigenvalue $\lambda_0(S)$ (given by $\nu_0(S)$ both described in point 5 of Proposition \ref{prop:descvaps} above), we can say that:
\begin{enumerate}
\item if $S<1$, we have $(u_0(x),v_0(x))=\left(\sin\left(\nu_0(1/2+x)\right),\sin\left(\nu_0(1/2-x)\right)\right)$, $x\in\left(-1/2,1/2\right)$;
\item if $S=1$, $(u_0(x),v_0(x))=(1+2x,1-2x),\ x\in\left(-1/2,1/2\right)$;
\item if $S>1$, $(u_0(x),v_0(x))=\left(\sinh\left(\nu_0'(1/2+x)\right),\sinh\left(\nu_0'(1/2-x)\right)\right)$, $x\in\left(-1/2,1/2\right)$,
where $\nu_0'=-i\nu_0\in\mathbb{R}$.
\end{enumerate}
This fact follows directly from point 2 of Proposition \ref{prop:vapsfups} and point 5 of Proposition \ref{prop:descvaps}. But we wanted to include this explicit description of $(u_0(x),v_0(x))$ here because it gives the asymptotic profile of the solutions (see Theorem \ref{thm:main}), and it can be interesting to have this explicit form if one is dealing with some applications of the problem.
\end{rem}

\begin{rem}
The sequence of $\nu_{n,j}$ can be naturally ordered by the increasing value of
the real part of its elements: $0\le \Re(\nu_{0})<\pi\le \ldots \leq \, n\pi<\Re(\nu_{n,1})\leq \Re(\nu_{n,2})<(n+1)\pi$ for all $n\ge 1$.
\end{rem}

\begin{rem}
As we will see in the proof of Proposition \ref{prop:descvaps}, the criterion for the subindex $j$ in $\nu_{n,j}$ and $\lambda_{n,j}$ will be the following one. We choose
$\Im(\nu_{n,j})>0$ if $j=1$, and $\Im(\nu_{n,j})<0$ if $j=2$. Then, we will see that the criterion happens to be the same for $\lambda_{n,j}$.
\end{rem}

\begin{rem}
The study of $\lambda_0(S)$ will be continued in Sections \ref{sec:domvap} and \ref{sec:decay}.  In Section \ref{sec:domvap} we will prove that it is indeed the dominant eigenvalue of the problem. More concretely, we will prove the inequality and gap condition $\Re(\lambda_{n,j}(S))<\lambda_0(S)-\varepsilon$ for some $\varepsilon=\varepsilon(S)>0$, and for all $n\geq 1$, $j=1,2$ (which is point 3 of Theorem \ref{thm:main}). Moreover, in Section \ref{sec:decay}, we will see that it gives the optimal decay rate of the solutions of problem \eqref{opA} (see the rest of Theorem \ref{thm:main}).
\end{rem}

\begin{figure}[htpb]
    \centering
   \includegraphics[width=0.4\textwidth]{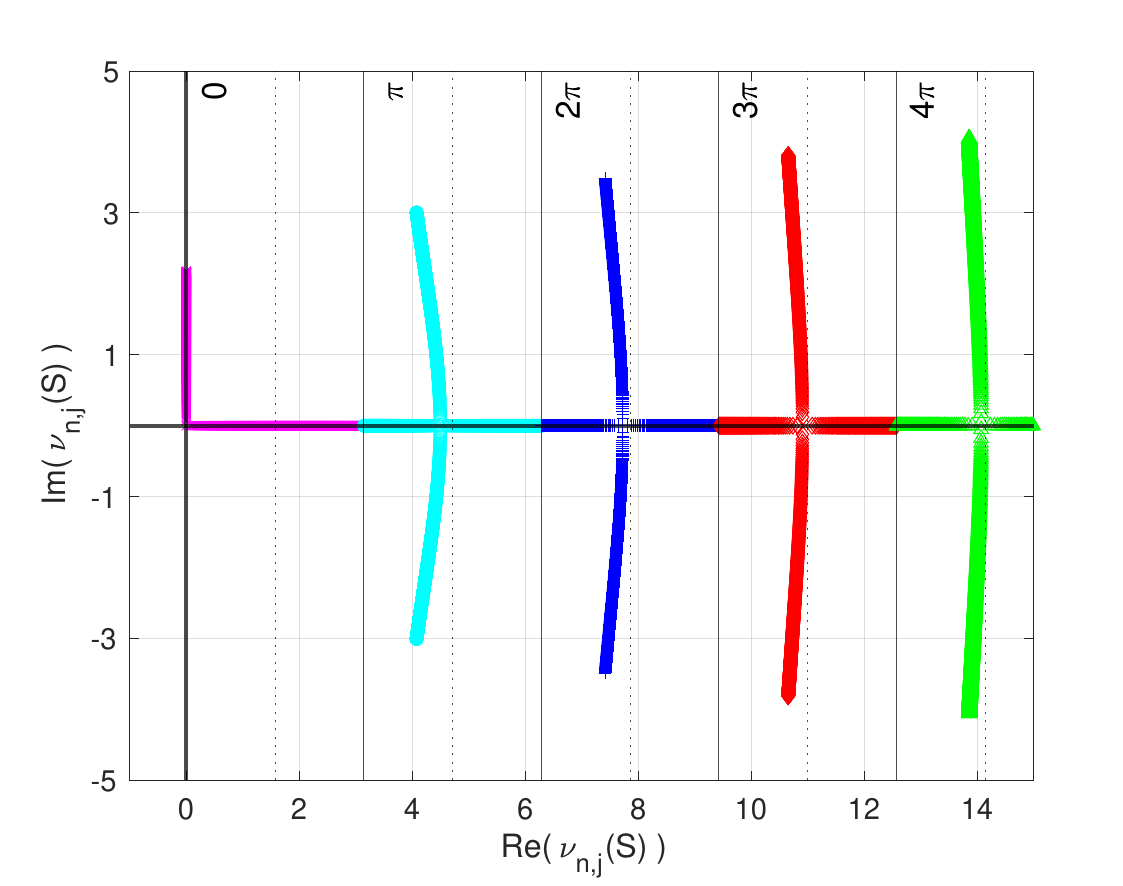}
   \includegraphics[width=0.4\textwidth]{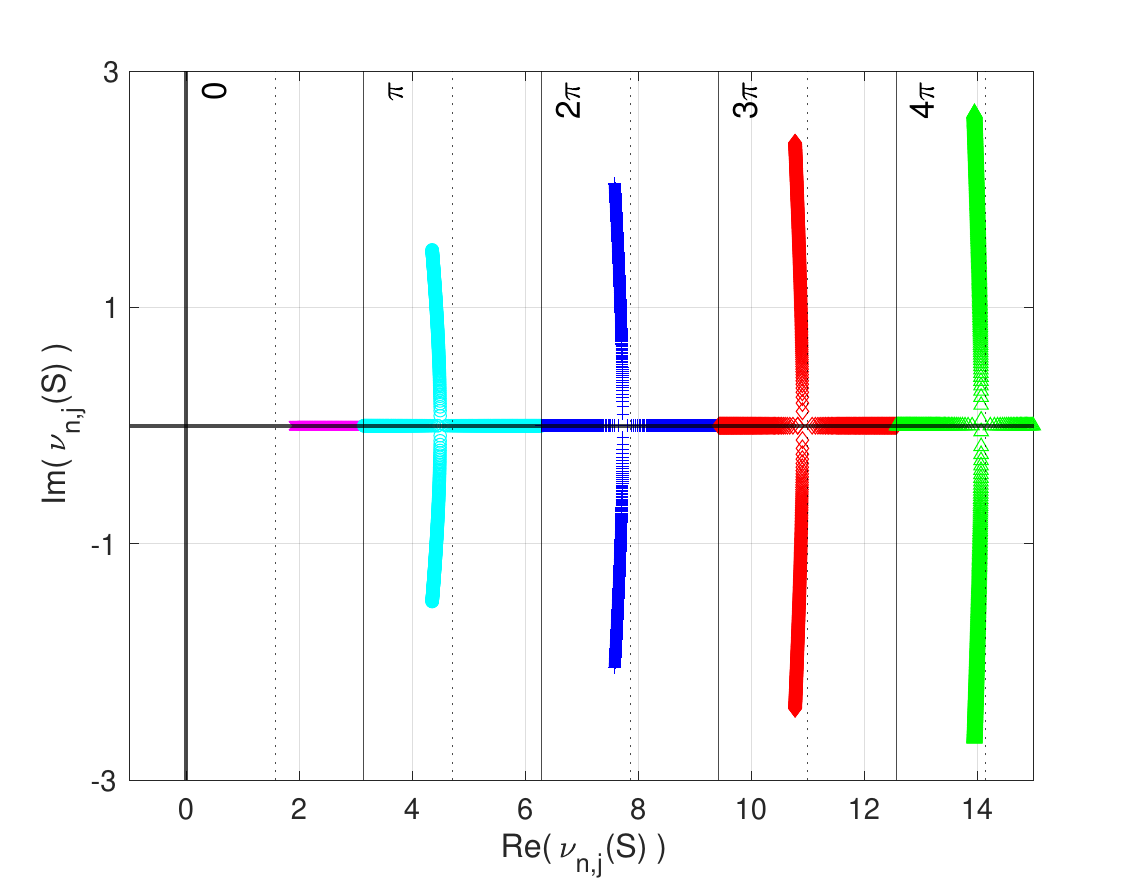}
      \includegraphics[width=0.4\textwidth]{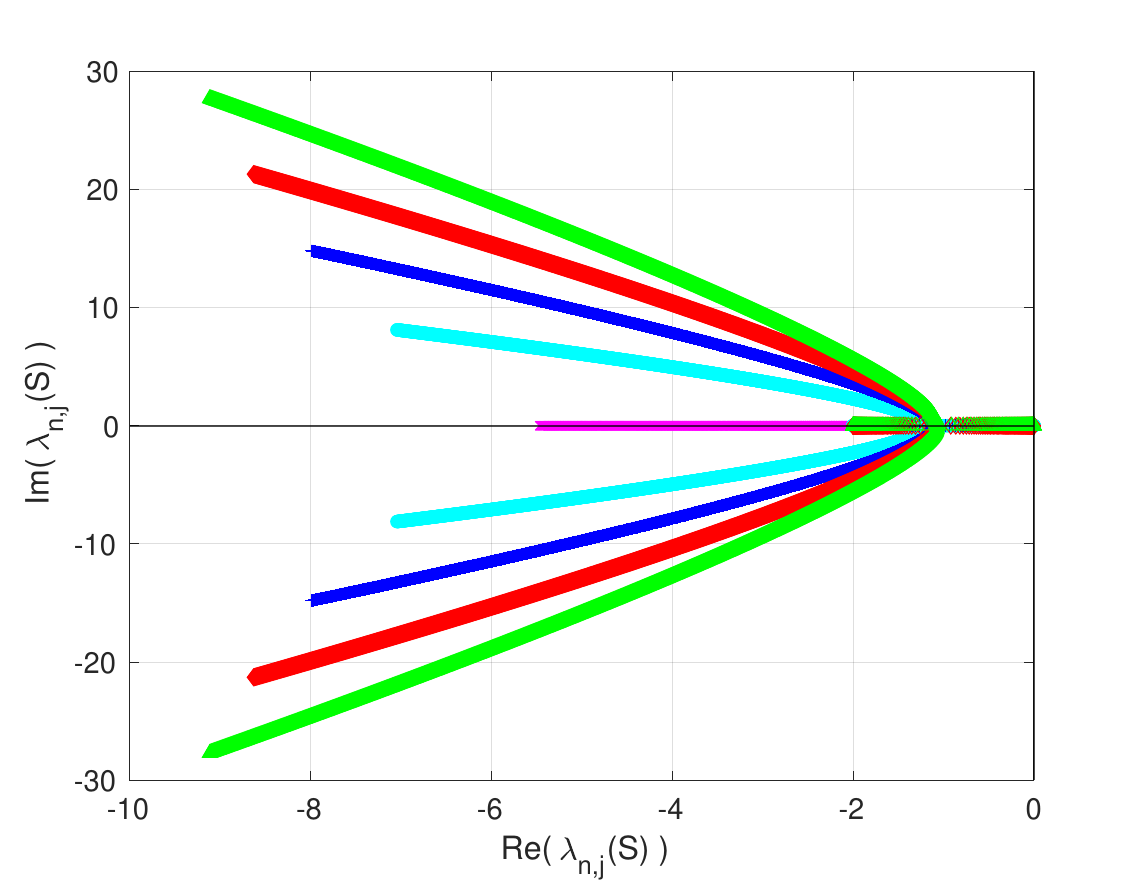}
      \includegraphics[width=0.4\textwidth]{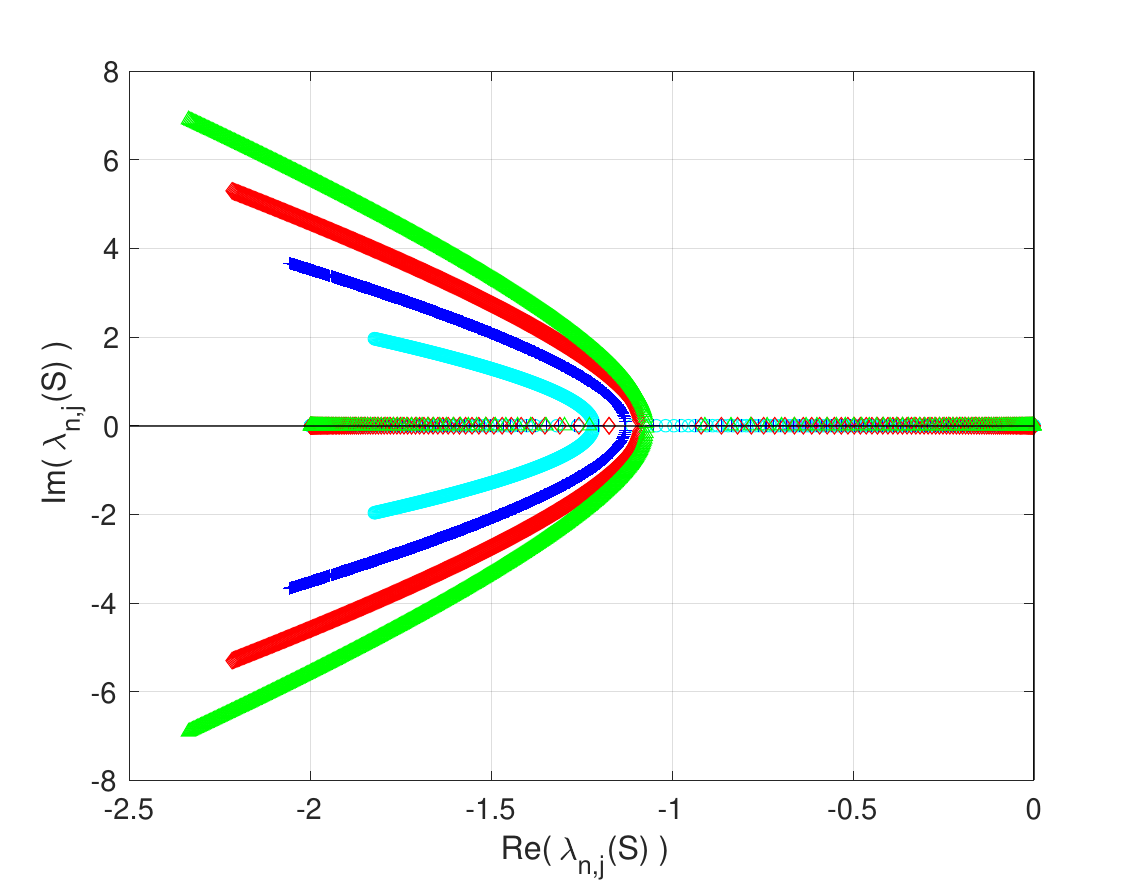}

    \caption{$\nu_{n,j}(S)$ (up) and $\lambda_{n,j}(S)$ (bottom) in the complex plane for $j=1,2,\ n=0,1,2,3,4$ (from left to right in the top row, and from the real axis up (or down), in the bottom row). We have $S\in(0,2]$ (first column) and $S\in(0,0.5]$ (second column). Note that the eigenvalues overlap on the real line. Also, as $S$ increases, the values of $\nu_{n,j}(S)$ and $\lambda_{n,j}(S)$ become non-real, and $\lambda_{n,j}(S)$ tend to $-\infty\pm i\infty$. }
    \label{fig:nu}
\end{figure}

\begin{figure}[htpb]
    \centering
   \includegraphics[width=0.4\textwidth]{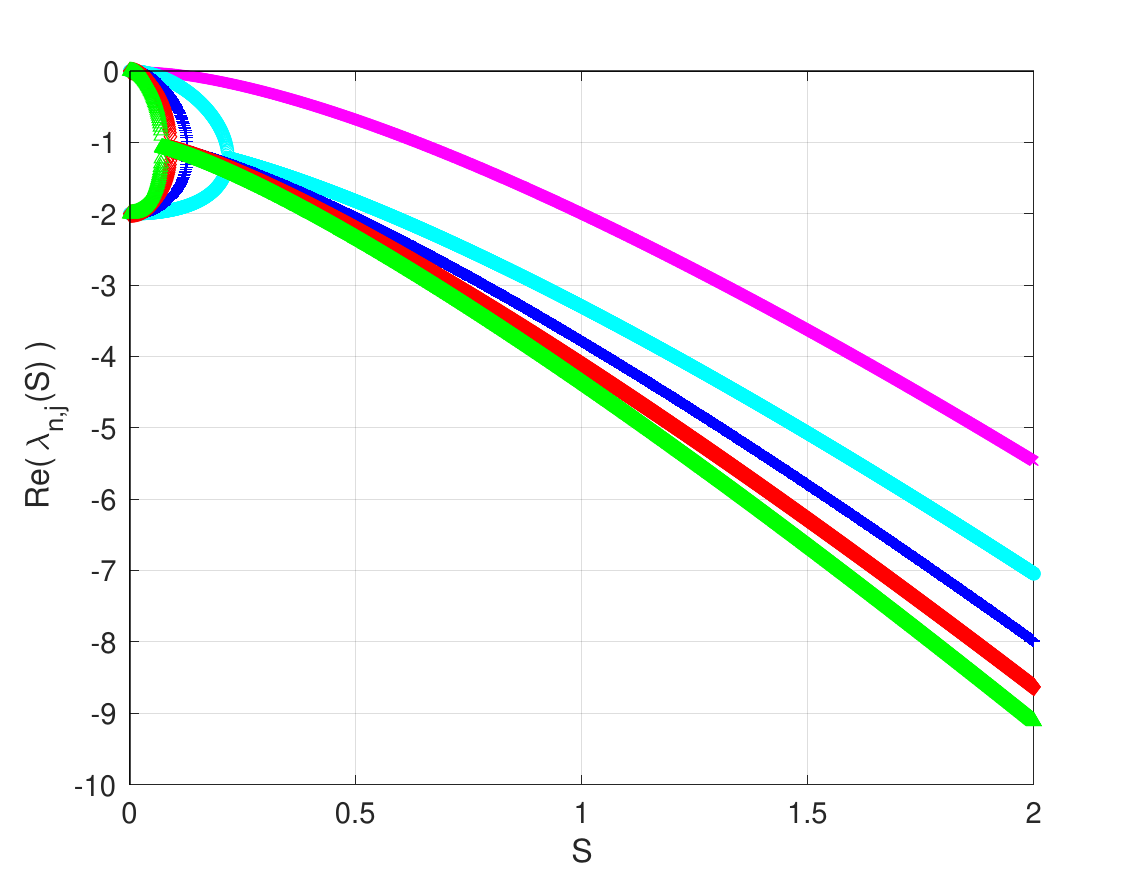}
   \includegraphics[width=0.4\textwidth]{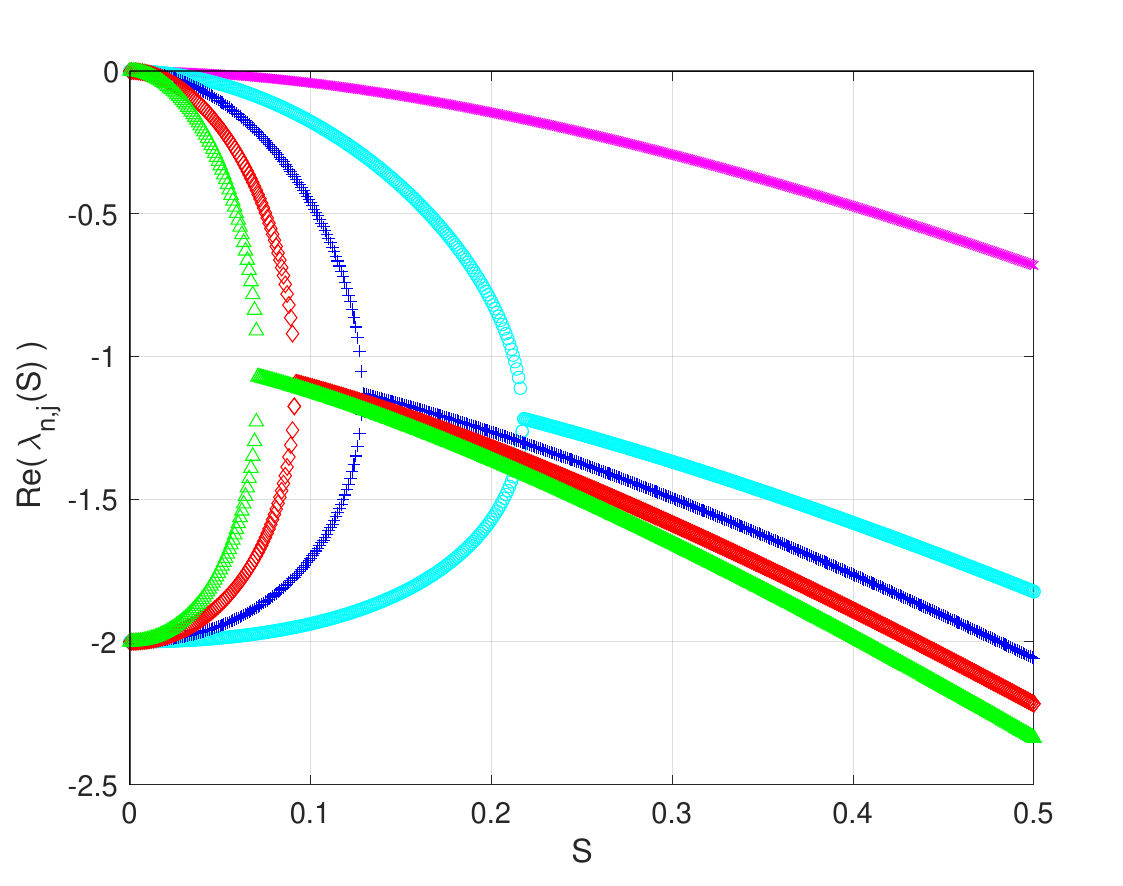}
      \includegraphics[width=0.4\textwidth]{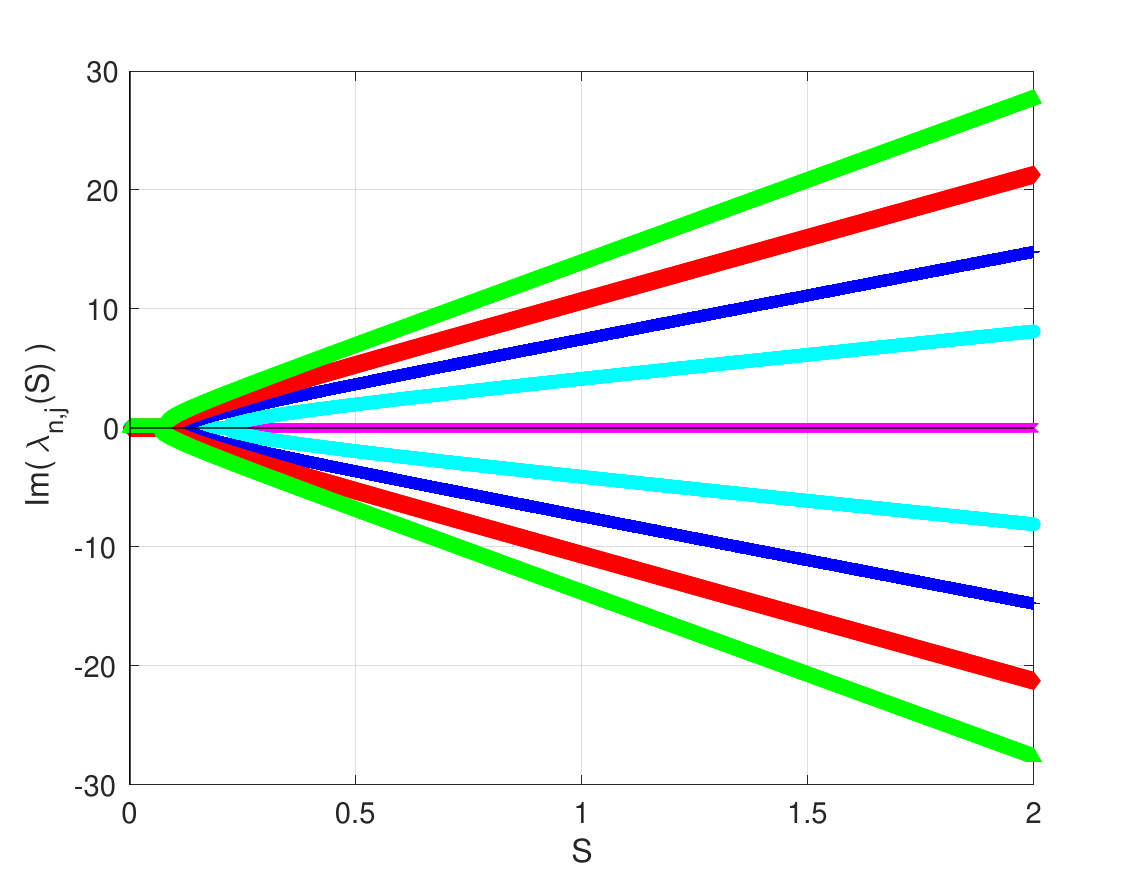}
   \includegraphics[width=0.4\textwidth]{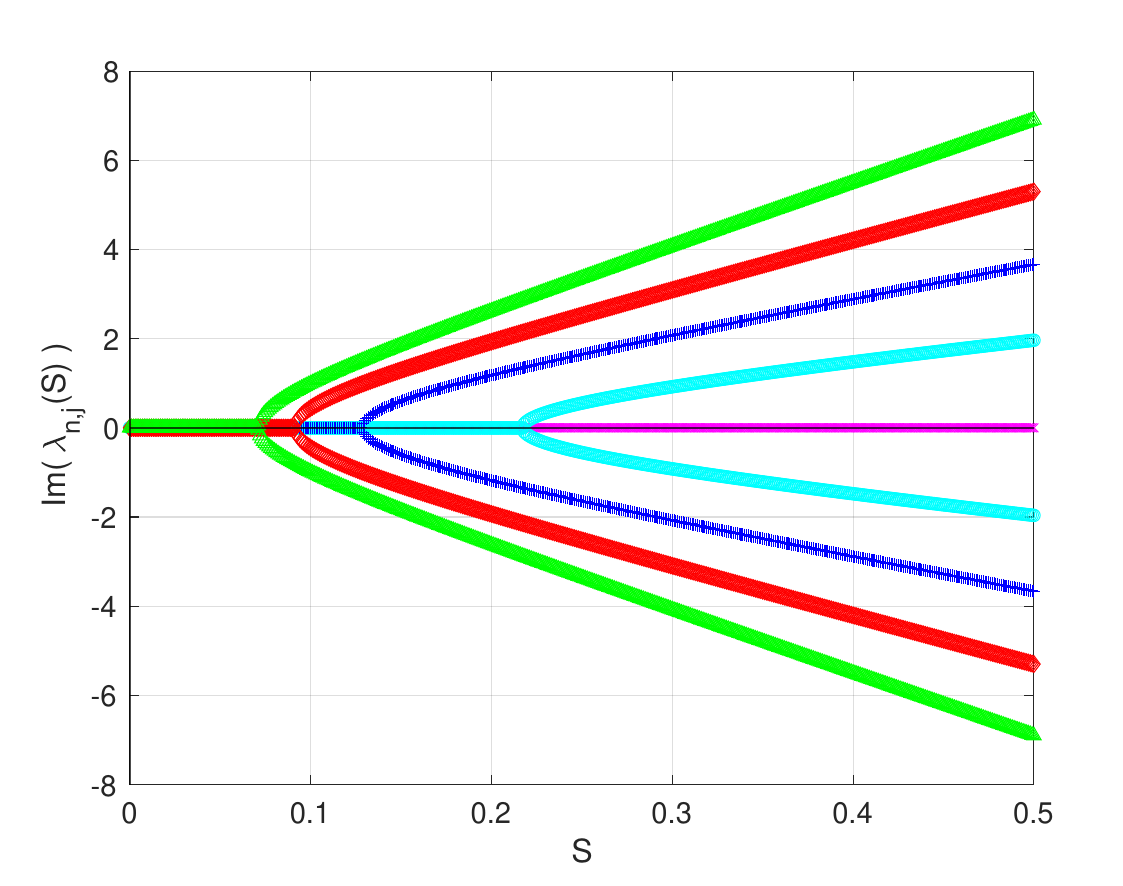}

    \caption{$\Re(\lambda_{n,j}(S))$ (top) and $\Im(\lambda_{n,j}(S))$ (bottom) for $j=1,2,\ n=0,1,2,3,4$ (top to down), with $S\in(0,2]$ (first column) and $S\in(0,0.5]$ (second column). Again, increasing $S$ makes the values of $\lambda_{n,j}(S)$, $n\geq 1$, become non-real and $\lambda_0(S) $ and $\Re(\lambda_{n,j}(S))$, $n\geq 1$, decrease to $-\infty$. In this sense, we say that $S$ is, indeed, a dissipative parameter: the larger it is, the faster the solutions tend to 0. }
    \label{fig:ReLambda}
\end{figure}

\begin{pf}
In Proposition \ref{prop:vapsfups} we have seen that the eigenvalues $\lambda\in\mathbb{C}\setminus\{-2\}$ of \eqref{opA} are the solutions of \eqref{ces} or \eqref{cea}, with $\nu=\sqr(\lambda^2-2\lambda)/S$, with $\nu\neq 0$. Looking at the graphs of the real functions $\sin\nu$ and $\pm S\nu$ for $S>0$, we can see that these solutions occur in pairs (see Figure \ref{fig:NUS1}). For small values of $S>0$ they are real and different, and belong to $(n\pi, (n+1)\pi)$ (with $n$ even in the case of a plus sign, and $n$ odd in case of a minus sign), may become equal for certain values of $S>0$, and become complex conjugate for larger values of $S>0$. As we have said in part 1 of the present proposition, we call them $\nu_{n,1}(S)$ and $\nu_{n,2}(S)$ ($n$ even in the symmetric case, $n$ odd in the antisymmetric one).

To be more precise, we call $S_m$ the value of $0<S < 1$ that makes either the equation $\sin(\nu)=S\nu$ or the equation $\sin(\nu)=-S\nu$ to admit a double real root $\nu$ in $[m\pi,(m+1/2)\pi]$, $m\geq 1$ (even or odd). Observe that the sequence of $S_m$ is decreasing in $m$, and that $S_m\to 0$ as $m\to\infty$. And that for all $S\geq S_{crit}=S_1$, all the solutions of $\sin(\nu)=S\nu$ or $\sin(\nu)=-S\nu$ are non-real  (except the solution $\nu_0(S)$ in $(0,\pi)$ for $S\in[S_{crit},1]$), with $S_{crit} \simeq 0.2172$ (see Figure \ref{fig:NUS1}). This proves parts 1 and 2 of this proposition.

We now consider the symmetric case, that is, the $\nu_{n,1},\nu_{n,2}$ ($n$ even) that are solution of the first equation of \eqref{ces}. Observe that for a fixed $S\in(S_{m+2},S_m]$ with $m\geq 2$ and even, we have $\nu_{n,1}(S),\nu_{n,2}(S)\in(n\pi,(n+1)\pi)$ for $1< n < m$, or $n=m$ and $S\neq S_m$, $n$ even, and that both are real and different. In fact, as $S$ increases, $\nu_{n,1}(S)$ increases, and $\nu_{n,2}(S)$ decreases, up to the value $S=S_m$, where $\nu_{m,1}(S)$ and $\nu_{m,2}(S)$ collide (this is in fact the double real root of $\sin(\nu)=S\nu$). In this last case, $\nu_{n,1}(S_n)=\nu_{n,2}(S_n)\in (n\pi, (n+1/2)\pi)$ ($n$ even).

\begin{figure}[htpb]
    \centering
    \includegraphics[angle=270, width=1.0\textwidth]{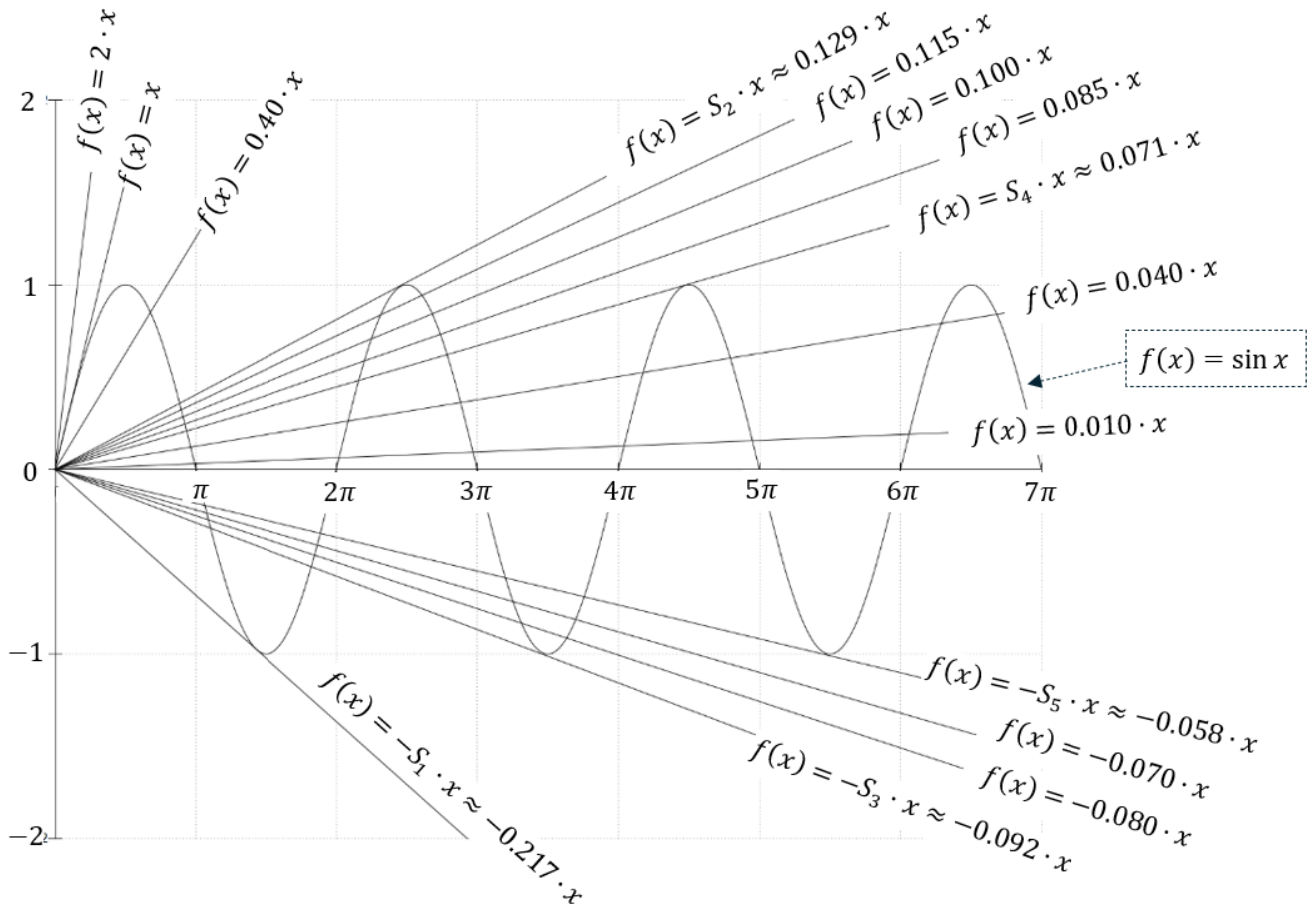}
    \caption{Intersections of the graphs of $\sin x$, and $Sx$ and $-Sx$ for different values of $S>0$}
    \label{fig:NUS1}
\end{figure}

Now, for a fixed $S\in(S_{m+2},S_m]$ with $m\geq 2$, and for $n>m$, the solutions $\nu_{n,1}(S)$ and $\nu_{n,2}(S)$ become non-real, the first with positive imaginary part, say, and the other its conjugate value.
In particular, if $S>S_{crit}=S_1$ this happens for all $n\geq 1$. These non-real $\nu_{n,1}$ and $\nu_{n,2}$ will necessarily be in the interior of the strip $Q_n$ defined in \eqref{eq:strip} and its conjugate $\overline{Q_n}$, respectively. That is,
\begin{equation*}
Q_n=\{z\in\mathbb{C}; \ n\pi<\Re(z)<(n+1/2)\pi, \ \Im(z)>0 \}.
\end{equation*}

We prove this fact by continuation on the parameter $S$, based on the so-called Argument Principle of Complex Analysis. One has to observe that along the boundary $\rm{Im}(z)=0$ of $Q_n$ the only possible solutions of $\sin\nu=S\nu$ are the $\nu_{n,1}(S)$ and $\nu_{n,2}(S)$, that exist only for $0<S\le S_m$ with $n\leq m$, already considered. Along the boundary $\Re(z)=n\pi$ of $Q_n$
one has $\sin z= i\sinh(\rm{Im}(z))$ (remember that $n$ is even), and a solution of $\sin\nu=S\nu$ would be impossible. Along the boundary $\Re(z)=(n+1/2)\pi$ of $Q_n$
one has $\sin z= \cosh(\rm{Im}(z))$, and a solution of $\sin\nu=S\nu$ would also be impossible. Finally, to see that the solutions $\nu_{n,1}$ cannot leave $Q_n$ through any upper boundary, and thus become unbounded, let us write $\nu=\alpha+i\beta$ where $\alpha$ and $\beta$ are real. Then, the equation $\sin(\nu)=S\nu$ can be written as the real system $\sin(\alpha)\cosh(\beta)=S\alpha$, and $\cos(\alpha)\sinh(\beta)=S\beta$. Squaring and adding the two equalities, and using that $|\sinh \beta|<\cosh \beta$ one gets the two inequalities $\sinh^2\beta<S^2(\alpha^2+\beta^2)<\cosh^2\beta$.

The first one is enough for our purposes since it means that $1<S^2((n+1/2)^2\pi^2+\beta^2)/\sinh^2\beta$, but this quotient tends to zero when $\beta\to\infty$. Therefore, for fixed $n$ it would be impossible to have a bounded sequence $S_k$, and the corresponding sequence $\nu_{n,1}(S_k)=\alpha_k+i\beta_k\in Q_n$ such that $\beta_k\to\infty$ when $k\to\infty$. The argument for $\nu_{n,2}$ and $\overline{Q_n}$ is the same.

Let us study now the antisymmetric eigenfunctions case, that is, $\nu$ satisfying the first equation of \eqref{eigenforms_asym},
$\sin\nu=-S\nu$.
Looking at the graphs of the real functions $\sin \nu$ and $-S\nu$ for $S>0$ (see Figure \ref{fig:NUS1}) observe that we obtain exactly the same results for $\nu_{n,j}(S)$ than for the symmetric case.

In order to prove \eqref{eq:nuasympt}, the asymptotic expansion of $\nu_{n,j}$ when $n\to\infty$ and $S>0$ fixed, we write $\nu_{n,j}=\alpha_{n,j}+i\beta_{n,j}$, with $\beta_{n,j}\neq 0$ (observe that when $n\to\infty$ this is always the case). The equation $\sin(\nu_{n,j})=\pm S\nu_{n,j}$ can be written as the real system of two equations
\begin{equation}\label{systemp}
\begin{cases}
\sin(\alpha_{n,j})\cosh(\beta_{n,j})=\pm S\alpha_{n,j}\\
\cos(\alpha_{n,j})\sinh(\beta_{n,j})=\pm S\beta_{n,j}.
\end{cases}
\end{equation}

Since $n\pi<\alpha_{n,j}<(n+1/2)\pi$ we have that $\alpha_{n,j}\to\infty$. As from the first equation of \eqref{systemp} we have  $\alpha_{n,j}=\pm\frac{1}{S}\sin(\alpha_{n,j})\cosh(\beta_{n,j})$, we obtain that $\beta_{n,j}\to\pm\infty$. In fact, recalling the sign criterion for the subindex $j$, we have $\beta_{n,1}\to\infty$ and $\beta_{n,2}\to-\infty$.

Then, necessarily $\lim_{n\to\infty} \frac{\beta_{n,j}}{\sinh(\beta_{n,j})}=0$, and using the second equation in the system \eqref{systemp}, we get $\lim_{n\to\infty}\cos(\alpha_{n,j})=0$. Since we know that $n\pi<\alpha_{n,j}<(n+1/2)\pi$, this necessarily implies that
$\alpha_{n,j}=(n+1/2)\pi+o(1)$ when $n\to\infty$.

Also, necessarily $\lim_{n\to\infty}|\sin(\alpha_{n,j})|=1$, or, in other words, $|\sin(\alpha_{n,j})|=1+o(1)$, and we can go to the first equation to obtain $\cosh(\beta_{n,j})=\frac{S\alpha_{n,j}}{|\sin(\alpha_{n,j})|}=S(n+1/2)\pi+o(1)$. Since $d/dy  \left({\rm argcosh}(y)\right)\to 0$ as $y\to\infty$, we deduce that
$\beta_{n,j}={\rm argcosh}(S(n+1/2)\pi)+o(1)=\pm\log\Bigl(S(n+1/2)\pi+\sqrt{S^2(n+1/2)^2\pi^2-1}\Bigr)+o(1)=\pm\log(2S(n+1/2)\pi)+o(1)=\pm\log(2Sn\pi)+o(1)$.
Again, recalling that we have chosen that $\Im(\nu_{n,1})>0$ and $\Im(\nu_{n,2})<0$, this completes the proof of \eqref{eq:nuasympt} and, hence, of part 3 of this proposition.

Once we have analysed $\nu_{n,j}(S)$ we can say something about the corresponding eigenvalues $\lambda_{n,j}$, given by the second equations of \eqref{ces} if $n$ is even, or \eqref{cea} if $n$ is odd. That is,  $\lambda_{n,j}=-1-\cos (\nu_{n,j})$, $j=1,2,\ n> 1$ even, or $\lambda_{n,j}=-1+\cos (\nu_{n,j})$, $j=1,2,\ n\geq 1$ odd. We can see that, if $\lambda_{n,j}\in\mathbb{C}\setminus\mathbb{R}$, the subindex $j$ in $\lambda_{n,j}\in\mathbb{C}\setminus\mathbb{R}$ actually follows the same criterion than in $\nu_{n,j}$: that is, $\Im(\lambda_{n,1})>0$ and $\Im(\lambda_{n,2})<0$ (to prove this, just recall that $\Re(\nu_{n,j}) \in (n\pi,(n+1/2)\pi)$ if $n\geq 1$).

When one knows $\sin(\nu)$ one can calculate $\cos(\nu)$ as $\pm\sqrt{1-\sin^2\nu}$, but it is still not clear which one of the two signs one has to take. In the case that $\nu_{n,1 }=\overline{\nu_{n,2}}\in\mathbb{C}\setminus\mathbb{R}$, we have that $n\pi<\rm{Re}(\nu)<(n+1/2)\pi$ since $\nu_{n,j}\in Q_n\cup \overline{Q_n}$. Therefore $\Re(\cos \nu)>0$ for $n$ even, and $\Re(\cos \nu)<0$ for $n$ odd. Then, if $n$ is even, we have $\cos(\nu_{n,j})=\sqr(1-\sin^2\nu_{n,j})$, and if $n$ is odd, we have $\cos(\nu_{n,j})=-\sqr(1-\sin^2\nu_{n,j})$. Using \eqref{ces} if $n$ is even, or \eqref{cea} if $n$ is odd, in both cases we obtain
$$\lambda_{n,j}(S)=-1-\sqr(1-S^2(\nu_{n,j}(S))^2),\ n\geq 1,\ j=1,2.$$

The above expression $\lambda_{n,j}(S)=-1-\sqr(1-S^2(\nu_{n,j}(S))^2)$ also holds if $S=S_m$ and $n=m$, as $\nu_{n,1}=\nu_{n,2}\in(n\pi,(n+1/2)\pi)$. But the same analysis cannot be applied exactly if $S\in(S_{m+2},S_m]$ and $n< m$, as the sign of $\Re(\cos \nu)= \cos \nu$ depends on the sub-indexes $n$ and $j$.
More concretely, since $\nu_{n,1}\in(n\pi,(n+1/2)\pi)$, the same kind of arguments above allow us to say that
$$\lambda_{n,1}(S)=-1- \sqr(1-S^2(\nu_{n,1}(S))^2), n\geq 1. $$
For $j=2$, however, the argument is a bit more tricky and the best we can say is
$$\lambda_{n,2}(S)=-1\pm \sqr(1-S^2(\nu_{n,2}(S))^2), n\geq 1.$$
Indeed, for most $S,n$, we have $\nu_{n,2}(S)\in((n+1/2)\pi,(n+1)\pi)$, hence, $\cos(\nu_{n,2}(S)) <0$ if $n$ even, and $\cos(\nu_{n,2}(S)) > 0$ if $n$ is  odd. That gives $\lambda_{n,2}(S)= -1+\sqr(1-S^2(\nu_{n,2}(S))^2) $ for these values of $S,n$.  But there are still some values of $S,n$ for which $\nu_{n,2}(S) \in (n\pi, (n+1/2)\pi)$ (essentially, those $S,n$ that make $\nu_{n,2}(S)$ almost equal to $\nu_{n,1}(S)$, that is $S$ close to $S_m$, $n$ close to $ m$). If this happens, $\cos\nu_{n,2}(S) > 0$ if $n$ even, and $\cos\nu_{n,2} (S)< 0$ if $n$ odd. This gives $\lambda_{n,2}(S)= -1-\sqr(1-S^2(\nu_{n,2}(S))^2) $ for those values of $S,n$. Note that when $S\to 0$, we have $\lambda_{n,2}(S)= -1+\sqr(1-S^2(\nu_{n,2}(S))^2) $ for a fixed $n\geq 1$. In particular, using the expressions above, it is direct to see that $\lambda_{n,1}(S)\to -2$ and $\lambda_{n,2}(S)\to 0$ when $S\to 0$.

\vspace{0.2cm}
We now proceed to prove inequality \eqref{lem:ineq}. If we express $\nu=x+iy$, the real part of the first equations of \eqref{ces} and \eqref{cea} can be expressed as $\sin(x) \cosh(y)=\pm  Sx$. Observe that if $\nu_{n,j}(S)\notin \mathbb{R}$ each of the previous equations can be written
as a function $y=y(x)$ with domain inside the strip $n\pi<x<(n+1)\pi$, with $n\geq 1$, and with $y>0$ (see point 3 of this proposition and Figure \ref{fig:yA} below). Observe  that  $y(x)\to\infty$ when $x\to n\pi$ or $(n+1)\pi$, and has a positive minimum at a certain point $e_n$.   Taking derivatives, we can see that $e_n=\tan(e_n)$, which ensures that
$e_n\in (n\pi, (n+1/2)\pi)$.

\begin{figure}[h]
    \centering
    \includegraphics[width=.9\textwidth]{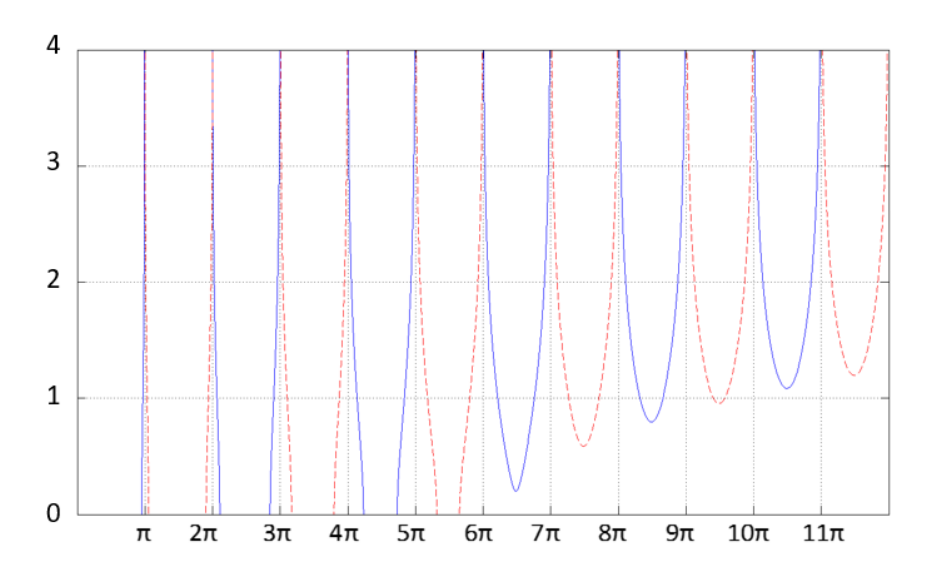}
    \caption{Graph of the function $y(x)$ for $S=0.05$. We observe that it is well defined and has a minimum in each $\left(n\pi, (n+1)\pi\right)$ from $n\geq 6$ (which corresponds to $\nu_{n,j}\notin\mathbb{R}$). The blue continuous line corresponds to the function obtained using \eqref{ces}, and the red discontinuous one corresponds to \eqref{cea}.}
    \label{fig:yA}
\end{figure}

Since the point $(x_{n,1},y_{n,1})$ has to be in the graph of $y(x)$, and using that $e_n$ is a minimum of $y(x)$ we have
$$\cosh(y_{n,1})\geq \cosh(y(e_{n}))=\pm\dfrac{Se_n}{\sin(e_n)}=\pm S\dfrac{e_n}{\tan(e_n)\cos(e_n)}=\pm \dfrac{S}{\cos(e_{n})}.$$
In addition, we know that $(x_{n,1},y_{n,1})$ is inside the strip $Q_n$ (see Proposition \ref{prop:descvaps}, point 3), which means that $x_{n,1}\in \left(n\pi,(n+1/2)\pi\right)$. So, $\cos(x_{n,1})<0$ if $n$ is odd, and $\cos(x_{n,1})>0$ if $n$ is even. This implies that $\pm \cos(x_{n,1})=|\cos(x_{n,1})|$ and thus $\Re (\lambda_{n,1})=-1\mp\cos(x_{n,1})\cosh(y_{n,1})=-1-|\cos(x_{n,1})|\cosh(y_{n,1})$. So we can write the last inequality as
$$\cosh(y_{n,1})\geq \dfrac{S}{|\cos(e_{n})|},$$
and therefore,
\begin{equation}\label{eq:re_bound2}
\Re(\lambda_{n,1})\leq-1-\frac{|\cos(x_{n,1})|}{|\cos(e_{n})|} S.
\end{equation}

We now want to prove that $x_{n,1}<e_n$ to continue inequality \eqref{eq:re_bound2}. To see this, notice that by dividing the first equation by the second one in \eqref{eq:system}  one has the \textsl{locus} of the intersections
$(x_{n,1},y_{n,1})$ as the curves
$$\dfrac{\tanh(y)}{y}=\dfrac{\tan(x)}{x}.$$
Observe that $\tanh(y)/y<1$ for all $y>0$. Hence, we need $\tan(x)/x<1$. Observe that this only holds for $n\pi < x < e_n$ (since $\tan(e_n)/e_n=1$, and $\tan(x)/x$ is an increasing function). Therefore, necessarily $x_{n,1}<e_n$. Then, from this, $|\cos(x_{n,1})|>|\cos(e_{n})|$, and we can state from \eqref{eq:re_bound2}
that inequality \eqref{lem:ineq} holds, that is,
\begin{equation*}\label{eq:max2a}
\Re(\lambda_{n,2}) =\Re(\lambda_{n,1})<-1-S,
\end{equation*}
for $\lambda_{n,j}\in\mathbb{C}\setminus\mathbb{R}$.

Regarding the asymptotic expressions of $\lambda_{n,j}$, we now observe that for large $S$ or $n$ the expansion of $\sqr(1-z^2)$ in terms of $1/z^2$ for $z$ large (where $z=S\nu_{n,j}$) allows us to say that \eqref{eq:lambdanu} is satisfied, that is
$$\lambda_{n,j}(S)= -1-(-1)^j S \nu_{n,j}\, i +o(1).$$
(again, recall the sign criterion of the imaginary part for choosing $j=1$ or $j=2$).

In particular, for $n\to\infty$, the leading order expansion of $\nu_{n,j}$ given in \eqref{eq:nuasympt} gives the leading order expansion of $\lambda_{n,j}$ given in \eqref{eq:lambdaasympt}, from where
$$\Re(\lambda_{n,j}(S)) \simeq -S\log(n) $$
and
$$\Im(\lambda_{n,j}(S)) \simeq -(-1)^j S n\pi  $$
(see Figure \ref{fig:ReLambda}).
This completes the proof of part 4 of this proposition.

\vspace{0.2cm}
The case $n=0$ is different. There is only one solution, namely $\nu_0(S)$, for all values of $S>0$, that denotes the solution of the symmetric case equation $\sin(\nu)=S\nu$ in the strip $\{ 0\leq \Re (z)\leq \pi, \ \Im(z)\geq 0\}$. Looking at the graph of $\sin(\nu)$ and $S\nu$, it is easy to see that $\nu_0(S)$ is real and belongs to $(0,\pi)$ if $0<S<1$, purely imaginary if $S>1$ (and we only consider the root with strictly positive imaginary part due to Remark \ref{rem:sqrt}), and equal to 0 if $S=1$. Indeed, in Proposition \ref{prop:vapsfups} we have seen that $\nu(S)=0$ does not have a non-zero eigenfunction unless $S=1$ and, hence, it is a spurious solution. In the particular case $S=\frac{2}{\pi}$ it is direct to see that $\nu_0(S)=\frac{\pi}{2}$.

If we follow this $\nu_0(S)$ for different values of $S$, we observe that we start with $\nu_0(S)\to \pi$ when $S\to 0$. Then it decreases until $\nu_0(S)=0$ when $S=1$ (actually with a triple collision of $\nu_0(1)=0$ with the spurious solution $\nu=0$ that we always have, and the also spurious $\nu=-\nu_0(1)$ coming from the left half plane). And, finally, for $S>1$, one still has the spurious solution $\nu=0$, but the other two solutions  $\nu_{0,1}(S)$ and $\nu_{0,2}(S)$ will be purely imaginary and of the form $\pm iy_0(S)$, where $\pm y_0(S)$ are the two real solutions of $\sinh y=Sy$, which exist only when $S>1$. Since we only consider those with positive imaginary part, we simply call it $\nu_0(S)$. Observe that $\Im(\nu_0(S))$ increases in $S$ if $S>1$.

Using the second equation of \eqref{ces}, we can see that $\lambda_0(S)=-1-\cos(\nu_0(S))$ is real and decreasing for all values of $S$. This is obvious for all $0<S\leq 1$. And for $S>1$, where $\nu_0(S)=i y_0(S)$ is purely imaginary, we use that $\cos(x+iy)=\cos x\cosh y-i\sin x\sinh y$. And, hence, $\lambda_0(S)=-1-\cosh(y_0(S))$ is also real and decreasing in $S$. The rest of the statements in part 5 of the proposition are also direct by using $\lambda_0(S)=-1-\cos(\nu_0(S))$.

This proves part 5 of the present proposition and, hence, completes the whole proof.

\end{pf}

\begin{prop}{\textbf {(Oscillation properties of the eigenfunctions)}}\label{prop:osc}

Let $U_{n,j}(x)=(u_{n,j}(x),v_{n,j}(x))$ be an eigenfunction of $A$ of eigenvalue $\lambda_{n,j}$. We claim that if $\lambda_{n,j}\in\R$, then $u_{n,j}(x)$ and $v_{n,j}(x)$ have exactly $n$ simple zeroes in the open interval $-1/2<x<1/2$. And if $\lambda_{n,j}\in\C\setminus\R$ then the two curves of the complex plane $z=u_{n,j}(x)$  and $z=v_{n,j}(x)$ satisfy that $d/dx\left(\arg (u_{n,j}(x))\right)<0$ and $d/dx\left(\arg(v_{n,j}(x))\right)>0$ for all $x\in(-1/2,1/2)$ if $j=1$ (the opposite if $j=2$). And
$\left|\lim_{x\to-1/2}\arg(u_{n,j}(x))-\arg(u_{n,j}(1/2))\right|$ and also $\left| \arg(v_{n,j}(-1/2))-\lim_{x\to 1/2}\arg(v_{n,j})\right|$ are both greater than $n\pi$ and less than $(n+1)\pi$, for $j=1,2$.

\end{prop}

\begin{figure}[h]
    \centering
    \includegraphics[width=0.9\textwidth]{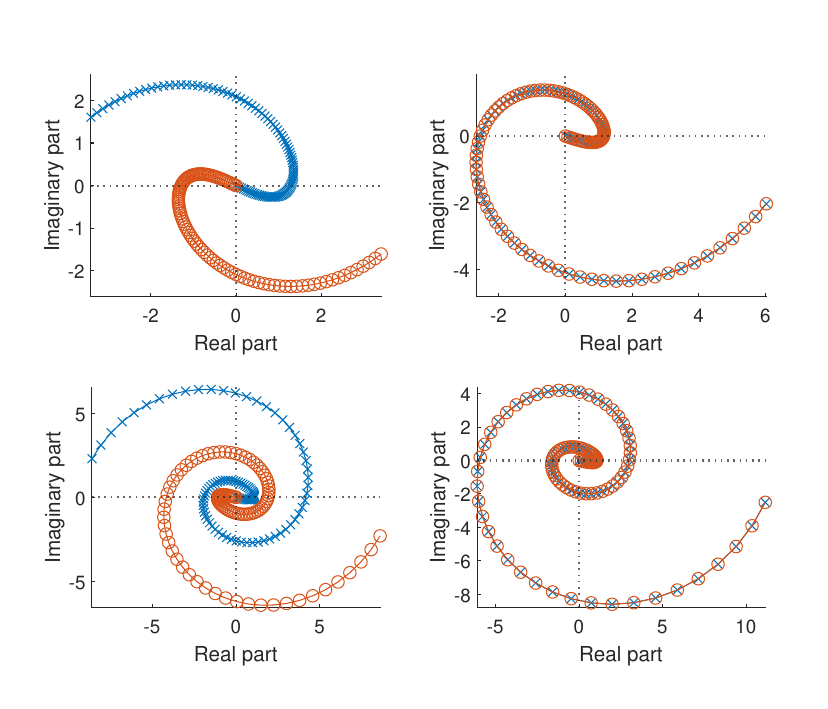}
    \caption{From left to right and from top to bottom, eigenfunctions $(u_{n,1},v_{n,1})$, $n=1,2,3,4$ for $S=0.8$ in the complex plane. In all cases, the blue cross line corresponds to $u_{n,1}$, and the red circle one to $v_{n,1}$ We can see that the $n$-th eigenfunction makes $n$ complete half-turns. The plots for the corresponding $(u_{n,2},v_{n,2})$ are similar.}
    \label{fig:oscil}
\end{figure}

\begin{rem}\label{rem:osc}
Note that the last statement of Proposition \ref{prop:osc} simply means that both curves $z=u_{n,j}(x)$ and $z=v_{n,j}(x)$ make $n$ complete half-turns around $z=0$ (see Figure \ref{fig:oscil}). This measure of the non-real oscillations, called rotation number by B. Simon in \cite{Simon2005}, also in the context of Sturm-Liouville Theory, is the natural generalization to the complex plane of the number of zeroes of real functions.

\end{rem}

\begin{pf}
All the eigenfunctions are complex scalar multiples of
$$(u_{n,j}(x),v_{n,j}(x))=\left(\sin(\nu_{n,j}(x-1/2)), (-1)^nu_{n,j}(-x)\right).$$
Since multiplication for a non-zero complex scalar is equivalent to a rotation and a homothety, it will suffice to prove our claim for $\sin(\nu_{n,j}(x-1/2))$. For simplicity, let us write $\nu_{n,j}=\alpha+i\beta$, $x-1/2=y$ and $u_{n,j}(x)=u(y)=\sin\left((\alpha+i\beta)y\right)$ with $0<y<1$.

If $\nu_{n,j}\in\R$ then $\beta=0$ and $n\pi<\alpha<(n+1)\pi$ (see Proposition \ref{prop:descvaps}), and $u(y)=\sin(\alpha y)$ has exactly $n$ simple zeroes in the open interval  $0<y<1$.

If $\nu_{n,j}\in\C\setminus \R$ then $u(y)=\sin(\alpha y)\cosh(\beta y)+i\cos(\alpha y)\sinh(\beta y)$. In particular, $|u(y)|^2>\sinh^2(\beta y)$, which cannot vanish if $y>0$. Also, the argument of $u(y)$, as of every complex number, is the imaginary part of its logarithm, and
$$
\dfrac{d}{dy}\arg(u(y))=\Im \left(\dfrac{u'(y)}{u(y)}\right)=\Im\left(\dfrac{u'(y)\overline{u(y)}}{|u(y)|^2}\right).$$
Then
$$\Im(u'(y)\overline{u(y)})=
-\cos(\alpha y)\sinh(\beta y)\bigl(\alpha\cos(\alpha y)\cosh(\beta y)+\beta\sin(\alpha y)\sinh(\beta y)\bigr)+$$
$$
\sin(\alpha y)\cosh(\beta y)\bigl(\beta\cos(\alpha y)\cosh(\beta y)-\alpha\sin(\alpha y)\sinh(\beta y)\bigr)=
$$
$$
-\alpha\sinh(\beta y)\cosh(\beta y)+\beta\sin(\alpha y)\cos(\alpha y)=
-\dfrac{1}{2}\alpha\sinh(2\beta y)+\dfrac{1}{2}\beta\sin(2\alpha y).
$$
We know that $\alpha>0$. Now suppose that also $\beta>0$ as well. Then, by using that $\sinh(2\beta y)>2\alpha\beta y$ (for $y>0$) and also $\beta\sin(2\alpha y)<2\beta\alpha y$ we conclude that the previous expression is strictly negative for all $y>0$. And strictly positive if $\beta<0$.

To see that the total increment of $\arg(u(y))$ for $0<y<1$ is greater than $n\pi$ and less than $(n+1)\pi$ it is sufficient to observe that
$$\Re(u(y))=\sin(\alpha y)\cosh(\beta y),$$ and that this real function has exactly $n$ simple zeroes in the open interval $0<y<1$, since $n\pi<\alpha<(n+1)\pi$. This concludes the proof.

\end{pf}

\section{The dominant eigenvalue $\lambda_0(S)$}\label{sec:domvap}
This section is devoted to prove point 3 of Theorem \ref{thm:main}. Recall that the number of real eigenvalues is finite for all the values of $S>0$, and that for the non-real ones we have the asymptotic behaviour for $n$ large given by formula \eqref{eq:lambdaasympt}. Hence, we can say that for each $R>0$ there are always only a finite number of eigenvalues with $\Re(\lambda(S))>-R$. Therefore, it will be enough to prove the strict inequality
\begin{equation}\label{eq:l0}
\Re(\lambda_{n,j}(S)) < \lambda_0(S) \hbox{  for all } n\ge 1, \, j=1,2,
\end{equation}
and that the geometric and algebraic multiplicity of $\lambda_0(S)$ is equal to 1 for all $S>0$. Both facts are needed to prove that $\lambda_0(S)$ is the dominant eigenvalue of our problem and it gives the optimal decay of the solutions, as Theorem \ref{thm:main} states and as we will prove in Section \ref{sec:decay}. For the proof, we use the results seen of Section \ref{sec:vapsfups}.

\begin{pf}\textbf{(of inequality \eqref{eq:l0})}

In the proof we distinguish three cases:

\begin{enumerate}
\item $0<S\le1$ and taking $\lambda_{n,j}(S)\in\mathbb{R}$, $n\geq 1,\ j=1,2$.

The existence of these real $\lambda_{n,j}(S)$ only happens if $S\le S_{crit}\simeq 0.2172$ and, in this case, $\nu_{n,j}(S)\in\mathbb{R}$ too (see Proposition \ref{prop:descvaps}). Then it is easy to see, just by looking at the graphs of $y=\sin(x)$ and $y=Sx$ that, in this case, $\nu_0\in (\pi/2,\pi)$ and, hence $\cos(\nu_0)<0$.

For positive real values of $\nu$, the first equations of \eqref{ces} and \eqref{cea} become  $S\nu=|\sin(\nu)|$.  Using this and that $S\nu_0<S\nu_{n,j}$ for $n\geq 1$ (see Proposition \ref{prop:descvaps}) one has also $|\sin(\nu_0)|<|\sin(\nu_{n,j})| $. Then,  $|\cos(\nu_0)|>|\cos(\nu_{n,j})|$, and finally
$$\lambda_{n,j}\le -1+|\cos(\nu_{n,j})|<-1+|\cos(\nu_0)|=-1-\cos(\nu_0)=\lambda_0.$$

\item $0<S\le 1$ and taking $\lambda_{n,j}(S)\in\mathbb{C}\setminus\mathbb{R}$, $n\geq 1,\ j=1,2$.

In this case, we have $\nu_{n,j}(S)\in\mathbb{C}\setminus\mathbb{R}$ too, and $\nu_0(S)\in\mathbb{R}$ (see Proposition \ref{prop:descvaps}). First, let us see that $\lambda_0(S) \geq -1-S$ when $0<S\leq 1$. Indeed, for $S=1$ the equality holds and is obvious as $\lambda_0(1)=-2$ (see Proposition \ref{prop:descvaps}, point 5). For $S<1$, we recall that $\nu_0$ is the real solution of $\sin(\nu)=S\nu$ in the interval $(0,\pi)$. There, we have $\sin(\nu_0)/\nu_0>\cos(\nu_0)$. Therefore, $\lambda_0= -1-\cos(\nu_0)>-1-\sin(\nu_0)/\nu_0= -1-S$.

This inequality combined with inequality\eqref{lem:ineq} in Proposition \ref{prop:descvaps} concludes that
$\Re(\lambda_{n,j})<\lambda_0$ also in the present case.

\item $S>1$.

Here, $\nu_0(S)$ and all the values of $\nu_{n,j}(S), \, n\geq 1,\, j=1,2$  are non-real (see Proposition \ref{prop:descvaps}). We write  $\nu=x+iy$ and the second equations of \eqref{ces} and \eqref{cea} become
\begin{equation}\label{eq:sys_vap}
\left\{ \begin{array}{lcl}
\Re(\lambda_{})=-1\mp \cos(x_{}) \cosh(y_{})\\
\Im(\lambda_{})=\pm\,\sin(x_{}) \sinh(y_{}),
\end{array}\right.
\end{equation}
which only need to be considered for $x\ge 0$ and $y>0$ (that is, $j=1$), since we are only interested in the real part of the eigenvalues. The values of $x,y$ must satisfy the first equations of \eqref{ces} and \eqref{cea}, which now become the system
\begin{equation}\label{eq:system}
\left\{ \begin{array}{c}
\sin(x) \cosh(y)=\pm Sx,\\
\cos(x)\sinh(y)=\pm Sy.
\end{array}\right.
\end{equation}

Since $ x_{0}=\Re(\nu_0)=0$, from the second equation in \eqref{eq:system} (with an absolute value) we have
\begin{equation}\label{eq:ineq}
\dfrac{\sinh(y_{0})}{y_{0}}=S<\dfrac{S}{|\cos(x_{n,1})|}
=\dfrac{\sinh(y_{n,1})}{y_{n,1}},\  \  n\geq 1,
\end{equation}
(remember that $x_{n,1}\in \left( \,n\pi,(n+1/2)\pi \,\right)$, see Proposition \ref{prop:descvaps}, point 3). And from the monotonicity of the function $\sinh(y)/y$ for $y>0$ we deduce that $y_{0}<y_{n,1}$, $n\geq 1$. Now, the monotonicity of the function $y/\tanh(y)$ gives
$$\dfrac{y_{0}}{\tanh(y_{0})}<\dfrac{y_{n,1}}{\tanh(y_{n,1})}\ \ \hbox{   and therefore  }
\ \ \cosh(y_{0}) < \dfrac{\sinh(y_{0})}{y_{0}}\dfrac{y_{n,1}}{\tanh(y_{n,1})},\ n\geq 1.$$

Now, by using the first equality of \eqref{eq:ineq}, the previous inequality can be written as
\begin{equation}\label{ineq2}
\cosh(y_{0}) < S\dfrac{y_{n,1}}{\sinh(y_{n,1})}\cosh(y_{n,1}).
\end{equation}
Finally, using \eqref{ineq2} and the last equality of \eqref{eq:ineq} in the first equation of \eqref{eq:sys_vap} (that can be written in the following form for $n$ odd and even), we have
$$\Re(\lambda_{n,2})= \Re(\lambda_{n,1})=-1-|\cos(x_{n,1})|\cosh(y_{n,1})<-1-\cosh(y_{0})=\lambda_0,$$
as we wanted to prove.

\end{enumerate}
\end{pf}

We now proceed to prove that the geometric and algebraic multiplicity of $\lambda_0(S)$ is equal to one. This is the second fact needed to prove the optimal decay of the solutions, as we will see in Section \ref{sec:decay}.

Roughly speaking, remember that the geometric multiplicity of an eigenvalue $\lambda$ of an operator $A$ is the dimension of $\textrm{Ker}(A-\lambda I)$. And that the algebraic multiplicity is the maximum $n\in\mathbb{N}$ such that $\textrm{dim}(\textrm{Ker}(A-\lambda I)^n) =\textrm{dim}(\textrm{Ker}(A-\lambda I)^{n+1})$. For operators with compact resolvent, this happens for $n$ large enough. See for exemple H. Brezis \cite{Brezis} for more details.
\begin{prop}\label{prop:multi}
\textbf{  }
\begin{enumerate}
\item The two components $(u_0,v_0)$ of the eigenfunction corresponding to $\lambda_0(S)$ can be chosen to be real and strictly positive for $-1/2<x<1/2$.
\item The multiplicity of $\lambda_0(S)$  (both geometric and algebraic) is equal to one.
\end{enumerate}
\end{prop}

\begin{pf}
\begin{enumerate}
\item From Proposition \ref{prop:vapsfups}, we know that one can take $u_0(x)=\sin(\nu_0(1/2+x))$, $x\in(-1/2,1/2)$ when $S\neq 1$. For $0<S<1$ we know that $0<\nu_0<\pi$ (Proposition \ref{prop:descvaps}, point 5), and, hence, $0<\nu_0(1/2+x)<\pi$. Hence, the result holds. For $S>1$ we know (again from point 5 of Proposition \ref{prop:descvaps}) that $\nu_{0}$ are purely imaginary, namely  $\nu_{0}=i\nu_0'$ with $\nu_0'>0$, and we can take $u_0(x)=\sinh(\nu_0'(1/2+x))$, that is strictly positive for $-1/2<x<1/2$.  For $S=1$ we know that $u_0(x)=1+2x$, $x\in(-1/2,1/2)$ (Proposition \ref{prop:vapsfups}), and the result also holds.The result also holds in all the cases for $v_0(x)$, since $v_0(x)=u_0(-x)$.

\item We know from part 2 of Proposition \ref{prop:vapsfups} that the geometric multiplicity of $\lambda_0(S)$ is equal to 1. To see that the algebraic multiplicity is also equal to 1, we will prove that it cannot be 2 or higher. That is, for the operator $A$ defined in \eqref{opA} we want to prove that, for the eigenfunction $(u_0,v_0)$ given before, the equation $(A-\lambda_0I)(\xi,\eta)^T=(u_0,v_0)^T$ has no solutions satisfying the boundary conditions $\xi(-1/2)=\eta(1/2)=0$. More explicitly,
\begin{equation}\label{multi}
\begin{cases}
-S\xi_x-(1+\lambda_0)\xi+\eta=u_0,\\
S\eta_x-(1+\lambda_0)\eta+\xi=v_0,\\
\xi(-1/2)=\eta(1/2)=0.
\end{cases}
\end{equation}

If we multiply both sides of the first equality by $v_0$ and both sides of the second by $u_0$, integrate by parts between $-1/2$ and $1/2$, using the boundary values, and we add the results, then we get
\begin{eqnarray*}
\int_{-\frac{1}{2}}^{\frac{1}{2}}\left[\xi(S(v_0)_x-(1+\lambda_0)v_0)+\eta v_0\right]dx+ \int_{-\frac{1}{2}}^{\frac{1}{2}}\left[\eta(-S(u_0)_x-(1+\lambda_0)u_0)+\xi u_0\right]dx = \\
= 2\int_{-\frac{1}{2}}^{\frac{1}{2}} u_0v_0\ dx.
\end{eqnarray*}
But since $(u_0,v_0)^T$ is an eigenfunction, this means that $S(v_0)_x-(1+\lambda_0)v_0+u_0=0$ and $-S(u_0)_x-(1+\lambda_0)u_0+v_0=0$, and the previous integral equality becomes
$$\int_{-\frac{1}{2}}^{\frac{1}{2}}\left[-\xi u_0+\eta v_0\right]dx+\int_{-\frac{1}{2}}^{\frac{1}{2}} \left[-\eta v_0+\xi u_0\right]=2\int_{-\frac{1}{2}}^{\frac{1}{2}} u_0v_0\ dx$$
which means $0=2\int u_0v_0\ dx$, and this is impossible because of part 1.

\end{enumerate}
\end{pf}
\begin{rem}
The previous proof of part 2 is reminiscent of the use of the adjoint operator in Fredholm's alternative. Both results given in parts 1 and 2 are in the spirit of the statement of the Krein-Rutman Theorem but in an explicit and simpler situation.
\end{rem}

\section{Optimal decay and asymptotic profile of the solutions}\label{sec:decay}

This section is of functional-analytic character, and is devoted to proving parts 1, 2, and 4 of Theorem \ref{thm:main}, and Theorem \ref{thm:compact}. That is the decay rate, its optimality, the asymptotic profile, and the compactness property of the solutions, respectively. The main tool will be the Spectral Mapping Theorem (SMT) in the following form.

\begin{prop}\label{prop:SMT}
We have
$$\sigma(T(t))=e^{t\sigma(A)}\cup\{0\},\ \textrm{ for } t>0$$
where $A$ is our operator, defined in \eqref{opA}, and $T(t), t\geq 0$, is the corresponding $C_0$-semigroup.
\end{prop}

The proof of this result will be given below. The SMT type results allow us to obtain information of $\sigma(T(t))$ (and, hence of the decay of the solutions) from the knowledge of $\sigma(A)$ (see \cite{EngelNagel}, for instance). Since we have described $\sigma(A)$ in Sections \ref{sec:vapsfups} and \ref{sec:domvap}, this kind of result seems particularly appropriate.

Before the SMT result, we need to begin with the proof of the compactness property of the semigroup given in Theorem \ref{thm:compact}. A basic step in this proof is the results of the classical paper of A.F. Neves, H. Ribeiro, and O. Lopes, \cite{Nevesetal86}. This is a well-known paper, and its results are still the basis of many recent works on hyperbolic systems, like, for example, J.E. Muñoz Rivera and M.G.Naso \cite{MRetal23}. Following this approach we obtain consequences essentially equal, but slightly better than those of M. Lichtner \cite{MLichtner2008}.

\begin{pf}(of Theorem \ref{thm:compact})

Our operator $A$ fits into the general form  of linear homogeneous hyperbolic systems given in \cite{Nevesetal86} formula (II.1) with, in their notation, $n=2$, $N=1$, $K=\begin{pmatrix}
S & 0\\
0 & -S
\end{pmatrix} $, $C=\begin{pmatrix}
1 & -1 \\
-1&1
\end{pmatrix}$, and $D=E=F=G=0$. According to this formula, an additional equation appears, namely $\dfrac{d}{dt}d(t)=0$ and the boundary condition $v(1/2,t)=d(t)$, which we can make compatible with our boundary conditions by taking $d(0)=0$ as the initial condition, as in the problem represented by formula (II.4) of the same reference. Following their reasoning, our problem can be compared to the simpler uncoupled system
\begin{equation}\label{simpler}
\begin{cases}
u_t+S u_x=-u\hspace{1cm} x\in(-1/2,1/2)\\
v_t-S v_x=-v\hspace{1cm} x\in(-1/2,1/2)\\
u(-1/2,t)=v(1/2,t)=0,
\end{cases}
\end{equation}
in the sense that if $T(t)$ is the semigroup generated by the original equation, and $T_0(t)$ is the semigroup generated by \eqref{simpler}, then $T(t)-T_0(t)$ is a compact operator for $t\ge 0$ (see Theorem A from \cite{Nevesetal86}).
But $T_0(t)$ can be explicitly calculated:
$$
T_0(t)\begin{pmatrix}
u_0(x)\\
v_0(x)
\end{pmatrix}=\begin{pmatrix}e^{-t} u_0(x-St) \hbox{ if } St-1/2<x, \hbox{ and $0$ otherwise}\\
e^{-t}v_0(x+St) \hbox{ if } St-1/2<-x, \hbox{ and $0$ otherwise}\end{pmatrix}$$
and we see that $T_0(t)=0$ for $t>1/S$. And the conclusion is that $T(t)$ is compact for $t>1/S$.
\end{pf}
We can now deduce the Spectral Mapping Theorem in the form stated above.
\begin{pf}\textbf{(of Proposition \ref{prop:SMT})}

We know from Theorem \ref{thm:compact} that $A$ generates an eventually compact semigroup. Therefore, Corollary IV 3.12 of \cite{EngelNagel} allows us to say that $\sigma(T(t))\setminus\{0\}=e^{t\sigma(A)}$, for $t\ge 0$.  Since we know that $\Re(\lambda_{n,j})\to -\infty$, $j=1,2$, when $n\to\infty$ (see Proposition \ref{prop:descvaps}, part 4b) we have that $0\in\sigma(T(t))$ for all $t>0$. This proves the SMT in the form given in the statement of this proposition.
\end{pf}

With all this, we are now ready to prove the optimal decay rate and the asymptotic profile of the solutions of our problem.

\begin{pf}\textbf{(of parts 1, 2, and 4 of Theorem \ref{thm:main})}

Again, as $A$ generates an eventually compact semigroup (see Theorem \ref{thm:compact}), Corollary IV 3.12 of \cite{EngelNagel} says that $s(A)=\omega_0(A)$, where $\omega_0(A)$ is the growth bound of $A$ (that gives the order of growth of $T(t)$), and $s(A)$ is the spectral bound of $A$. We know from Lemma \ref{lem:spec} that $A$ has only point spectrum. This implies that the spectral bound is $s(A)=\lambda_0(S)$, the dominant eigenvalue of $A$ described in Sections \ref{sec:vapsfups} and \ref{sec:domvap}. As $\omega_0(A)=s(A)$, this means that for each $\varepsilon>0$ there exists an $M\ge1$ such that
$$\|T(t)\|\le Me^{t(\lambda_0(S)+\varepsilon)} \textrm{ for all }t\ge 0.$$

But we can now improve this decay estimate, since we know that $\lambda_0(S)$ is geometrically and algebraically simple and that $\sup\left\{\Re\left\{\sigma(A)\setminus\{\lambda_0\}\right\}\right\}<\lambda_0$ (see Lemma \ref{lem:spec}, inequality \eqref{eq:l0} and Proposition \ref{prop:multi}) .

To see this, we can divide the set $\sigma(A)$ by a vertical line $\rm{Re}(z)=\omega$ in such a way that $\lambda_0>\omega$ and $\rm{Re}(\lambda)<\omega$ for all $\lambda\in\sigma(A)\setminus\{\lambda_0\}$, for some $\omega>0$. Let us call $\Pi_0$ and $I-\Pi_0$ the spectral projections associated with this decomposition of $\sigma(A)$. Since these projections commute with $A$ we have that $T(t)=e^{tA}=e^{t\Pi_0(A)}+e^{t(I-\Pi_0)A}$. The first semigroup is one-dimensional and is simply given by $e^{t\Pi_0(A)}U_0=e^{t\lambda_0}\Pi_0(U_0)$. Applying \cite{EngelNagel}, Corollary IV 3.12, again to $e^{t(I-\Pi_0)A}$ and using $\sup\{\rm{Re}\{\sigma(A)\setminus\{\lambda_0\}\}\}<\omega$ we see that

$$\|T(t)U_0\|=\| e^{t\Pi_0(A)}U_0+e^{t(I-\Pi_0)A}U_0\| \le M\left(
e^{\lambda_0 t}\|\Pi_0 U_0\| + e^{\omega t}\|(I-\Pi_0)U_0\|\right) \le M'
e^{\lambda_0 t}\|\Pi_0 U_0\|.$$

From this we conclude that $\|T(t)\|\le M'e^{\lambda_0(S)\, t}$, as stated in point 1 of Theorem \ref{thm:main}. The optimality of this inequality (point 2 of Theorem \ref{thm:main}) comes from the fact the  equality holds for solutions whose initial condition is a multiple of the eigenfunction associated to $\lambda_0(S)$.

From the previous inequality we also deduce point 4 of Theorem \ref{thm:main}: under the generic assumption $\Pi_0U_0\ne 0$, that is satisfied by an open and dense subset of $X$, one has $T(t)U_0=e^{\lambda_0\,t}\Pi_0 U_0+O(e^{t\omega})$. This shows that these solutions approach a multiple of the first eigenfunction $(u_0(x),v_0(x))$ faster (at a rate $e^{\omega t}$ or even faster) than the rate at which they approach zero (that is, $e^{t\lambda_0}$) (remember that $\omega<\lambda_0<0$). Therefore, the first eigenfunction $(u_0(x),v_0(x))$ can be properly called the asymptotic profile of the solutions.

\end{pf}

\section{Acknowledgements}
J. Menacho is part of the Catalan research group 2021 SGR 01228 recognized by AGAUR.
M. Pellicer and J. Solà-Morales are supported by the Spanish grant PID2021-123903NB-I00 funded by MCIN/AEI/10.13039/501100011033 and by ERDF ``A way of making Europe'', by the Spanish grant RED2022-134784-T funded by MCIN/AEI/10.13039/501100011033, and by the Catalan grant 2021-SGR-00087.

\end{document}